\documentclass[final,a4paper,leqno]{siamltex}
\usepackage[english]{babel}
\usepackage{amsmath,amssymb,amsfonts,graphicx,subcaption,outlines}
\usepackage{myshortcuts,float}
\usepackage{xcolor,bm}
\pagecolor{white}
\usepackage{fancyhdr}
\usepackage{MnSymbol}
\usepackage{csvsimple}
\usepackage{tikz}
\usepackage{pgfplots}
\usetikzlibrary{decorations.pathreplacing}

\DeclareMathOperator*{\argmin}{argmin}
\DeclarePairedDelimiter\abs{\lvert}{\rvert}%
\newcommand\vectorize{\mathord{\mathrm{vec}}}
\usepackage{appendix}
\usepackage{comment}
\newcommand{\norm}[1]{\left\lVert#1\right\rVert}
\usepackage[linesnumbered,lined,boxed,commentsnumbered,ruled,vlined]{algorithm2e}
\let\oldnl\nl
\newcommand{\nonl}{\renewcommand{\nl}{\let\nl\oldnl}}

\usepackage[hidelinks]{hyperref}
\usepackage[labelfont=bf]{caption}
\newcommand*\samethanks[1][\value{footnote}]{\footnotemark[#1]}
\iftrue 

  \usepackage{tikz}
  \usepackage{pgfplots}
  \pgfplotsset{
    compat=newest,
    tick label style={font=\scriptsize},
    label style={font=\scriptsize},
    legend style={font=\scriptsize}
  }
   \usetikzlibrary{decorations.pathreplacing}
  \iffalse 
     \usepgfplotslibrary{external} 
     
  \else
     \renewcommand{\tikzsetnextfilename}[1]{}
  \fi
\else
  \usepackage{tikzexternal}
  \tikzsetexternalprefix{gfx_tmp/}
  
\fi

\newtheorem{remark}[theorem]{Remark}

\title{Preconditioned infinite GMRES for parameterized linear systems}

\author{Siobhán Correnty\thanks{Department of Mathematics, Royal Institute of Technology (KTH), SeRC Swedish
e-Science Research Center, Lindstedtsvägen 25, Stockholm, Sweden
\{\texttt{correnty,eliasj}\}\texttt{@kth.se}} \and
Elias Jarlebring\samethanks \and Kirk M. Soodhalter\thanks{School of Mathematics, Trinity College Dublin, The University of Dublin, College Green, Dublin 2, Ireland \texttt{ksoodha@maths.tcd.ie}}}

\begin{document}
\captionsetup[subfigure]{justification=centering}

\include{title}
\maketitle
\begin{abstract}
We are interested in obtaining approximate solutions to parameterized linear systems of the form $A(\mu) x(\mu) = b$ for many values of the parameter $\mu$. Here $A(\mu)$ is large, sparse, and nonsingular, with a nonlinear analytic dependence on $\mu$. Our approach is based on a companion linearization for parameterized linear systems.
The companion matrix is similar to the operator in the infinite Arnoldi method [E. Jarlebring, W. Michiels, K. Meerbergen, \textit{Numer. Math.}, 122(1):169–195, 2012], and we use this to adapt the flexible GMRES setting [Y. Saad, \textit{SIAM J. Sci. Comput.}, 14(2):461–469, 1993]. In this way, our method returns a function $\tilde{x}(\mu)$ which is cheap to evaluate for different $\mu$, and the preconditioner is applied only approximately. This novel approach leads to increased freedom to carry out the action of the operation inexactly, which provides performance improvement over the method infinite GMRES proposed in [Jarlebring, Correnty, accepted for publication in \textit{SIAM J. Matrix Anal. Appl.}, 2022], without a loss of accuracy in general. We show that the error of our method is estimated based on the magnitude of the parameter $\mu$, the inexactness of the preconditioning, and the spectrum of the linear companion matrix. Numerical examples from a finite element discretization of a Helmholtz equation with a parameterized material coefficient illustrate the competitiveness of our approach. The simulations are reproducible and publicly available online.
\end{abstract}

\begin{keywords}
inexact preconditioning, parameterized linear systems, Krylov methods, companion linearization, shifted linear systems, infinite Arnoldi, inner stopping criteria
\end{keywords}

\begin{AMS}
15A06, 65F08, 65F10, 65F50, 65N22, 65P99
\end{AMS}

\section{Introduction}
We develop techniques for solving linear systems where the system matrix and the solution depend on a parameter $\mu$, i.e., systems of the form
\begin{align} \label{eq:our-prob}
	A(\mu) x(\mu) = b,
\end{align}
where $\mu \in \mathbb{R}$, $A(\mu) \in \mathbb{R}^{n \times n}$ nonsingular, analytic, and nonlinear in $\mu$, and $b \in \mathbb{R}^n$. In particular, we are interested in the solution for many different values of $\mu$. There are many applications which can be described by this general setting. The matrix $A(\mu)$ can, for example, stem from a discretization of a differential operator, a setting arising in the study of partial differential equations (PDEs) with uncertainty. See Section~\ref{sec:Helmholtz-sec} for an example involving a finite element discretization of Helmholtz equation with parameterized material coefficient $\mu$. Linear systems of the form \eqref{eq:our-prob} also occur in the study of dynamical systems, in particular in the computation of the transfer function, where $\mu$ is the Laplace variable. See Section~\ref{sec:time-delay} for an example of this arising from a time-delay system.

One approach for solving \eqref{eq:our-prob} for many values of $\mu$ could be to solve the systems corresponding to each $\mu$ of interest separately. However, as we increase the accuracy of the discretization, the dimension of the systems increases. Consequently, the cost of solving each system increases at least linearly with the dimension, even with state of the art iterative methods for sparse matrices. The method proposed in this paper is based on the well-established Krylov subspace method GMRES \cite{Saad:1986:GMRES} along with its flexible variant \cite{Saad1993}, combined with algorithms from nonlinear eigenvalue problems (NEPs). We summarize the main idea of our method here.

Under the assumptions above, the matrix $A(\mu)$ can be expressed locally by an infinite Taylor series expansion centered around origin, i.e.,
\begin{align} \label{eq:A-taylor}
	A(\mu) = \sum_{\ell=0}^{\infty} A_{\ell} \mu^{\ell},
\end{align}
where
\begin{align} \label{eq:taylor-coeff}
	A_{\ell} := A^{(\ell)}(0) \frac{1}{\ell !} \in \mathbb{R}^{n \times n}
\end{align}
are the Taylor series coefficients. In our setting, we assume further that the Taylor coefficients in \eqref{eq:A-taylor} do not vanish after a certain degree, and many of the derivatives of $A(\mu)$ are computationally available. This situation arises when, for example, $A(\mu)$ can be expressed as a sum of products of matrices and functions $A(\mu)=C_1f_1(\mu)+\ldots+ C_kf_k(\mu)$, where $k \ll n$. The method proposed here efficiently approximates the solution to \eqref{eq:our-prob} for many values of the parameter $\mu$ simultaneously, without truncation in the Taylor series expansion \eqref{eq:A-taylor}. Our novel approach offers significant improvements in performance over the prior work \cite{JarlebringCorrenty1} for the considered examples.


The nonlinear dependence on the parameter $\mu$ in \eqref{eq:our-prob} is addressed with a technique called companion linearization, commonly used in the study of polynomial eigenvalue problems (PEPs), see e.g. \cite{Mackey:2006:VECT}, but has also been used for paramerized linear systems \cite{GuSimoncini}. It can be summarized as follows. We denote the truncated version of \eqref{eq:our-prob} as
\begin{align} \label{eq:our-prob-trunc}
	A_m (\mu) x_m (\mu) = b,
\end{align}
where
\begin{align} \label{eq:A-trunc}
	A_m(\mu) := \sum_{\ell = 0}^m A_{\ell} \mu^{\ell}
\end{align}
is assumed to be nonsingular and $A_{\ell}$ as in \eqref{eq:taylor-coeff}, and consider the system of equations
\begin{align} \label{eq:lin}
	\left(
	\begin{bmatrix}
	A_0 & A_1 & A_2 & \cdots & A_m \\
	& I & & & \\
	& & I & & \\
	& & &\ddots & \\
	& & & & I
	\end{bmatrix}
	- \mu
	\begin{bmatrix}
	0 & 0& \cdots & \cdots & 0 \\
	I & & & & \\
	& I & & & \\
	& & \ddots & & \\
	& & & I & 0
	\end{bmatrix}
	\right)
	\begin{bmatrix}
	x_m(\mu) \\
	\mu x_m(\mu) \\
	\mu^2 x_m(\mu) \\
	\vdots \\
	\mu^m x_m(\mu)
	\end{bmatrix}
	=
	\begin{bmatrix}
	b \\
	0 \\
	0 \\
	\vdots \\
	0
	\end{bmatrix}.
\end{align}
Equation \eqref{eq:lin} is of the form
\begin{align} \label{eq:unshifted-trunc}
	\left(  K_m - \mu M_m \right) v_m (\mu) = c_m,
\end{align}
where $K_m$, $M_m \in \mathbb{R}^{(m+1)n \times (m+1)n}$ are constant matrices and $c_m \in \mathbb{R}^{(m+1)\times n}$ is a constant vector. Here the matrix $K_m$ contains the coefficients of the truncated Taylor series expansion $A_m(\mu)$ as in \eqref{eq:A-trunc}. The solutions to \eqref{eq:our-prob-trunc} and the solutions to \eqref{eq:lin} are equivalent in a certain sense, which we prove in Section~\ref{sec:comp-lin}. Note that the parameter $\mu$ appears only linearly in \eqref{eq:unshifted-trunc}.

The block triangular matrix $K_m$ has inverse $K_m^{-1}$, expressed exactly as
\begin{align} \label{eq:Km-inv}
	K_m^{-1} =
	\begin{bmatrix}
	A_0^{-1} & -A_0^{-1} A_1 & \cdots & -A_0^{-1} A_m  \\
	& I & & \\
	& & \ddots & \\
	& & & I
	\end{bmatrix} \in \mathbb{R}^{(m+1)n \times (m+1)n},
\end{align}
and we consider an equivalent formulation of \eqref{eq:unshifted-trunc} with right preconditioning, i.e.,
\begin{align} \label{eq:shifted-trunc}
	\left( I - \mu M_m K_m^{-1} \right) \hat{v}_m(\mu) = c_m,
\end{align}
where $\hat{v}_m(\mu) = K_m v_m(\mu)$. The matrix $K_m^{-1}$ in \eqref{eq:shifted-trunc} is referred to as the preconditioning matrix.

The linear system in \eqref{eq:shifted-trunc} incorporates a shift with the identity matrix. Approaches for solving systems of the same form as \eqref{eq:shifted-trunc} have been developed in previous works, which we now briefly summarize. See, e.g., \cite{ShiftedFreund, FrommerMaass99:20,GuSimoncini,KilmerOleary01,etna_vol45_pp499-523}, as well as \cite{Bahkos:2017:MULTIPREC} with multiple preconditioners, \cite{Baumann2015NestedKM} with a nested framework, and \cite{FrommerGlassner98} with restarting. As Krylov subspaces are invariant under scaling of the matrix as well as addition of a scalar multiple of the identity, applying the method GMRES to the system in \eqref{eq:shifted-trunc} requires the construction of an orthogonal basis matrix for the Krylov subspace given by
\begin{align} \label{eq:our-krylov}
	\mathcal{K}_k \left(M_m K_m^{-1}, c_m \right) =\text{span} \{ c_m, M_m K_m^{-1} c_m, \ldots, (M_m K_m^{-1})^{k-1} c_m \}.
\end{align}

The basis matrix for the Krylov subspace in \eqref{eq:our-krylov} is built using the infinite Arnoldi method (see, for example, \cite{Jarlebring:2012:INFARNOLDI,TensorArnoldi}). This method constructs the basis independently of the truncation parameter $m$. In this way, our method approximates solutions to \eqref{eq:lin} with $m\to\infty$, i.e., without neglecting any terms in the Taylor series expansion of $A(\mu)$. Equivalently, the proposed method approximates solutions to \eqref{eq:our-prob}.

We build the basis matrix for the Krylov subspace in \eqref{eq:our-krylov} once. After constructing this basis matrix with the infinite Arnoldi method, approximating the solution to \eqref{eq:lin}, and equivalently \eqref{eq:our-prob}, reduces to solving one small least squares problem for every value of $\mu$ considered. In this way, our method returns a parameterization of the solution to \eqref{eq:our-prob}.

Additionally, this work proposes a method which applies the preconditioner inexactly within a flexible framework. Here the action of $A_0^{-1}$ in $K_m^{-1}$ can be applied approximately when constructing the orthogonal basis matrix of \eqref{eq:our-krylov}. Moreover, the accuracy of the application can be relaxed according to a dynamic inner stopping criteria, inspired by \cite{inexKry}. In general, the preconditioner is applied almost exactly when the residual of the outer method is large and with decreasing accuracy as the residual is reduced. In practice, the level of accuracy can be relaxed dramatically without degrading convergence.

A flexible preconditioning setting has been studied in \cite{Notay2000,vanDerVorstVuik1994,vuik1993} as well as in \cite{ElmanErn2001} using multigrid as a preconditioner, in \cite{WarsaBenzi} with 2-stage preconditioning and in \cite{BaglamaCal1998,Erhel1996,Joubert1994}, where the application of preconditioner was modified based on information from previous iterations. Additionally, in \cite{Chapman1998,SimonciniSzyld2003}, where Krylov methods were used as preconditioners and in \cite{Axelsson1991,Szyld2001}, where stopping criteria was utilized. The theory for inexact Krylov methods has also been investigated in several prior works. Specifically, in \cite{Golub1999InexactPC,inexKry}, as well as in the series of papers \cite{Bouras2000ARS,Bouras00arelaxation}, matrix-vector multiplications within Krylov methods were performed in a way such that the inexactness of the products increased as the methods progressed, without disrupting convergence.

The main contribution of this paper is a new algorithm which we call inexact infinite GMRES. This method offers considerable improvement over the previously presented infinite GMRES \cite{JarlebringCorrenty1}. The advances are specifically with respect to the most expensive part of the algorithm, the application of $A_0^{-1}$. Here the flexible framework allows for a cheap, inexact application of the action of $A_0^{-1}$, as this is a component of the preconditioning matrix. Numerical experiments show that inexact infinite GMRES has substantially lower complexity than infinite GMRES, without a loss of accuracy. Thus, the newly proposed algorithm solves \eqref{eq:our-prob} for many different values of $\mu$ with greater efficiency.

While inexact infinite GMRES uses the infinite Arnoldi method to build the Krylov basis matrix, there are many algorithms for NEPs that could have been considered, for instance, the compact rational Krylov method \cite{VanBeeumen:2015:CORK} or the two-level orthogonal Arnoldi method \cite{Kressner:2014:TOAR}. Additionally, successful approaches for solving the Helmholtz equation include preconditioning combined with fast iterative solvers, such as in \cite{Erlangga08,ErlanggaEtAl2004}, or through multiple iterative solves, as in \cite{BaylissEtAl1983,ElmanErn2001}.

Alternative strategies have been proposed within the framework of model order reduction.
Such approaches were developed in prior works, such as \cite{SturlerEtAl15,FengEtAl2013}, as well as in \cite{AntoulasEtAl2020} and references therein and in \cite{AhujaEtAl2012} incorporating recycling. In \cite{KressnerToblerLR}, low-rank tensor decompositions were utilized to solve parameterized linear systems. As these methods represent a completely different approach than the ones proposed in this paper, we have chosen not to make a direct comparison to them. Additionally, the technique of Krylov recycling has been investigated in other works on solving linear systems in general (see \cite{KilmerSturler06,ParksEtAl2006,SOODHALTER2014105}).

The paper is organized as follows. Section~\ref{sec:comp-lin} shows the linearization used in the method, along with a proof of the equivalence of \eqref{eq:our-prob-trunc} and \eqref{eq:unshifted-trunc}. In Section~\ref{sec:algorithm} we summarize the standard flexible GMRES method and present the algorithm inexact infinite GMRES. The dynamic inner stopping criteria is also explained in this section, as well a proof showing the inexact application of the preconditioner does not disrupt the convergence of the new method. Section~\ref{sec:convergence-the} presents convergence theory for the proposed method based on the standard theory for the method GMRES. Section~\ref{sec:simulations} illustrates the theory with an example from time-delay systems. This section also highlights the competitiveness of our approach, shown by an example with a finite element discretization of a Helmholtz equation.

In this work, we denote the vectorization of a matrix $W = \begin{bmatrix} w_0, w_1, \ldots, w_m \end{bmatrix} \in \mathbb{R}^{n \times (m+1)}$ by
\begin{align} \label{eq:vectorize}
	\vectorize(W) := \begin{bmatrix} w_0^T, w_1^T, \ldots, w_m^T \end{bmatrix}^T \in \mathbb{R}^{n(m+1)},
\end{align}
and we define the family of matrices indexed by $m \in \mathbb{Z}_{+}$
\begin{align} \label{eq:tilde-D-m}
	\tilde{D}_m :=
	\begin{bmatrix}
	0 & 1 &  & & \\
	   &    & \ddots & & \\
	   & & &1 & 0
	\end{bmatrix} \in \RR^{(m-1) \times (m+1)},
\end{align}
where $\tilde{D}_1 := 0 \in \RR$. The $m$th largest singular value of a matrix $A \in \mathbb{R}^{n_1 \times n_2}$ is given by $\sigma_m(A)$.

\section{Companion linearization} \label{sec:comp-lin}
Companion linearization has been used in prior works on PEPs (see, e.g., \cite{Mackey:2006:VECT}), and in works for solving linear systems \cite{GuSimoncini,JarlebringCorrenty1}. Linearizing a system of the form \eqref{eq:our-prob} allows us to exploit properties of Krylov subspaces, leading to computationally competitive approaches. More precisely, we build the basis matrix for the Krylov subspace in \eqref{eq:our-krylov} once and reuse it to solve the parameterized system for many values of $\mu$.

Furthermore, the companion linearization in \eqref{eq:lin} has a certain structure, which allows us to apply the infinite Arnoldi method when building the basis matrix. The result of this choice is an algorithm independent of the truncation parameter $m$, explained completely in Section~\ref{section:inf-flex-gmres}. The following theorem shows an equivalence in solutions of the truncated parameterized system \eqref{eq:our-prob-trunc} and the linearized equation system \eqref{eq:lin}.

\begin{theorem} \label{theorem:linearization}
Let $A(\mu)$ be as in \eqref{eq:our-prob} with truncated expansion as in \eqref{eq:A-trunc}. Then, the solutions to \eqref{eq:our-prob-trunc} and \eqref{eq:unshifted-trunc} are equivalent in the following sense.
\begin{enumerate}
\item[(a)] If solutions to \eqref{eq:unshifted-trunc} are of the form
\begin{align} \label{eq:our-vm}
	v_m(\mu) = \begin{bmatrix} \mu^0 x_m^T(\mu), \mu x_m^T(\mu), \ldots, \mu^m x_m^T(\mu)  \end{bmatrix}^T \in \mathbb{R}^{(m+1)n},
\end{align}
then $x_m(\mu)$ satisfies \eqref{eq:our-prob-trunc}
\item[(b)] If $x_m(\mu)$ satisfies \eqref{eq:our-prob-trunc}, then $v_m(\mu)$ as in \eqref{eq:our-vm} satisfies \eqref{eq:unshifted-trunc}.
\end{enumerate}
\end{theorem}

{\em Proof.}
\begin{enumerate}
\item[(a)] Assume \eqref{eq:unshifted-trunc} holds. Then, the first block row in \eqref{eq:unshifted-trunc} implies \eqref{eq:our-prob-trunc} is fulfilled.
\item[(b)] Assume \eqref{eq:our-prob-trunc} holds. Then, first block row in \eqref{eq:unshifted-trunc} satisfied. The $\ell$th block row gives
\[
	\mu^{\ell-1} x_m(\mu) - \mu \mu^{\ell-2} x_m(\mu) = 0
\]
for $\ell=2,\ldots,m+1$. \hfill $\endproof$
\end{enumerate}
\begin{remark}[Structure of the companion linearization]
The companion linearization presented in \eqref{eq:lin} has a different structure than the linearizations used in prior works on the infinite Arnoldi method. The structure here was specifically developed in order to incorporate right preconditioning in a flexible setting as proposed in \cite{Saad1993}, while preserving the growing nonzero structure of the Krylov vectors. Additionally, the linearization in \eqref{eq:lin} allows for a dramatic reduction in complexity by applying the action of $A_0^{-1}$ inexactly while still taking advantage of the shift invariance property of Krylov subspaces, as can be seen in the equivalent preconditioned formulation in \eqref{eq:shifted-trunc}. This is described completely in Section~\ref{sec:algorithm}. The structure of the companion linearization in prior works did not meet this criteria.
\end{remark}
\begin{remark}[Truncation parameter $m$]
Note that taking $m \rightarrow \infty$ implies that \eqref{eq:our-prob} and \eqref{eq:lin} are equivalent in an analogous sense. In practice, Algorithm~\ref{alg:FIGMRES} performs computations with the same number of Taylor series coefficients as the desired dimension of the Krylov subspace, described completely in Section~\ref{sec:Algder}.
\end{remark}
\begin{remark}[Unique solution to preconditioned system]
Note that Theorem~\ref{theorem:linearization} implies the system matrix in \eqref{eq:unshifted-trunc} is nonsingular. Equation \eqref{eq:shifted-trunc} is obtained by multiplication with matrices $K_m^{-1}$ and $K_m$, nonsingular by assumptions on $A(\mu)$. Thus, the transformed system matrix in \eqref{eq:shifted-trunc} is also nonsingular.
\end{remark}

\section{Algorithm} \label{sec:algorithm}

\subsection{Flexible GMRES} \label{sec:FGMRES}
As a preparation for the derivation of the proposed algorithm, we consider flexible preconditioning of the GMRES algorithm \cite{Saad1993}. For better efficiency of the linear system solver, Krylov subspace methods are commonly carried out on a preconditioned system which exhibits better convergence characteristics than the original system. This requires additional solutions to inner systems, which if selected well, results in a more efficient method. This is the foundation for the method flexible GMRES (FGMRES), which solves the resulting inner systems inexactly, thus resulting in faster convergence. FGMRES is summarized below. For a full explanation of the standard GMRES algorithm, see \cite{Saad:1986:GMRES}.

At the $j$th iteration of the standard GMRES algorithm performed on the system $Ax = b$, with $A \in \mathbb{R}^{n \times n}$, $b \in \mathbb{R}^n$, we form the Arnoldi factorization
\[
	AQ_j = Q_{j+1} \underline{H}_j,
\]
where $Q_j = \begin{bmatrix} q_1,\ldots,q_j \end{bmatrix} \in \mathbb{R}^{n \times j}$ is a matrix whose columns form an orthonormal basis for the Krylov subspace of dimension $j$ associated with the matrix $A$ and vector $b$ (see \eqref{eq:our-krylov}) and $\underline{H}_j \in \mathbb{R}^{(j+1) \times j}$ is upper Hessenberg. In practice, we perform one matrix-vector product, $y = A q_j$ with $q_1 = b/ \norm{b}$, and orthogonalize this vector against the columns of $Q_j$ by a Gram--Schmidt process. This new vector is used to form $Q_{j+1}$. The orthogonalization coefficients are stored in $\underline{H}_j$. We compute the $j$th iterate $x_j$ as
\[
	x_j = x_0 + Q_j w_j,
\]
where
\begin{subequations}
\begin{eqnarray} \label{eq:comp-y}
	w_j &=& \argmin_{w \in \mathbb{R}^j} \norm{\beta e_1 - \underline{H}_j w}_2, \\
	\beta &=& \norm{b-Ax_0}_2,
\end{eqnarray}
\end{subequations}
and $x_0$ is the initial iterate.
The basic idea behind right preconditioned Krylov subspace methods is to solve a modified system of the form
\[
	A M^{-1} (Mx) = b,
\]
where $M \in \mathbb{R}^{n \times n}$ nonsingular. In general, we do not need to form the matrix $AM^{-1}$ explicitly. The action of this matrix can instead be computed by solving the linear system
\begin{align} \label{eq:pre-con-step}
	Mz_i=q_i
\end{align}
for $i=1,\ldots,j$ and on iteration $j$ forming the Arnoldi factorization
\begin{align} \label{eq:pre-arn}
	A M^{-1} Q_j = Q_{j+1} \underline{H}_j.
\end{align}
The $j$th iterate of preconditioned GMRES is thus given by
\[
	x_j = x_0 + M^{-1} Q_j w_j,
\]
where $w_j$ is as in \eqref{eq:comp-y} and $\underline{H}_j$ as in \eqref{eq:pre-arn}. For efficiency, it is important that \eqref{eq:pre-con-step} is easy to compute for each vector $q_i$, for example by means of an iterative method.

The FGMRES method instead considers changing the preconditioner on each step of the Arnoldi algorithm, i.e., \eqref{eq:pre-con-step} is replaced by
\begin{align} \label{eq:pre-con-step2}
	M_i \tilde{z}_i=q_i
\end{align}
for $i=1,\ldots,j$. On iteration $j$ we form the Arnoldi factorization
\begin{align} \label{eq:flex-arn}
	A \tilde{Z}_j = Q_{j+1} \underline{H}_j,
\end{align}
where $\tilde{Z}_j = \begin{bmatrix} \tilde{z}_1,\ldots, \tilde{z}_j \end{bmatrix}$ and $q_1 = b/ \norm{b}$. Thus, the $j$th iterate of FGMRES is given by
\[
	x_j = x_0 + \tilde{Z}_j w_j,
\]
where $w_j$ is computed as in \eqref{eq:comp-y} and $\underline{H}_j$ is as in \eqref{eq:flex-arn}.

In this paper, we consider the solution to a system shifted by the identity matrix and scaled with $-\mu$. Therefore, we form the shifted Arnoldi factorization
\begin{align*}
	(I - \mu A M^{-1}) Q_j = Q_{j+1} (\underline{I}_j - \mu \underline{H}_j)
\end{align*}
in the general case and
\begin{align} \label{eq:shifted-arn}
	Q_j - \mu A \tilde{Z}_j = Q_{j+1} (\underline{I}_j - \mu \underline{H}_j)
\end{align}
in the flexible form, where
\begin{align} \label{eq:Im}
	\underline{I}_j \in \RR^{(j+1) \times j}
\end{align}
is an identity matrix with an extra row of zeros. Here we have used the scaling- and shift-invariance properties of Krylov subspaces, i.e.,
\begin{align} \label{eq:shift-inv}
	\mathcal{K}_k (AM^{-1},b) = \mathcal{K}_k (I - \mu A M^{-1}, b).
\end{align}
The $j$th iterate of FGMRES for the shifted system is thus given by $x_j$, where
\begin{subequations} \label{eq:der-least-squares}
\begin{eqnarray}
	w_j &=& \argmin_{w\in\mathbb{R}^j} \norm{\beta e_1 - (\underline{I}_j - \mu \underline{H}_j)w}_2, \\
	x_j &=& x_0 + \tilde{Z}_j w_j,
\end{eqnarray}
\end{subequations}
and $\underline{H}_j$ is as in \eqref{eq:shifted-arn}. Thus, the flexible method requires the extra storage of the matrix $\tilde{Z}_j \in \mathbb{R}^{n \times j}$ but no additional cost in computation compared to preconditioned GMRES.

\subsection{Inexact infinite GMRES} \label{section:inf-flex-gmres}
\subsubsection{Algorithm derivation} \label{sec:Algder}
We specialize the flexible GMRES algorithm described in Section~\ref{sec:FGMRES} for shifted systems to the system in \eqref{eq:unshifted-trunc} with an inexact application of the preconditioner $K_m^{-1}$. The structure of our linearization allows for an efficient computation of the corresponding Arnoldi factorization. Here we show the action of an operator of size $(m+1)n \times  (m+1)n$ can be achieved using vectors and matrices of size $n \times n$ and smaller.

The standard preconditioned GMRES algorithm applied to \eqref{eq:unshifted-trunc} requires the Arnoldi factorization
\begin{align} \label{eq:shift-arn-our}
	\left( I - \mu M_m K_m^{-1} \right) Q_j = Q_{j+1} \left( \underline{I}_j - \mu \underline{H}_j \right),
\end{align}
formed by computing the product
\begin{align} \label{eq:multMK}
	M_m K_m^{-1} q_i
\end{align}
on the $i$th iteration of GMRES. This step is followed by an orthogonalization procedure and an update of the approximation via the computation of a least squares solution analogous to \eqref{eq:der-least-squares}. Here $q_1 = c_m/\norm{c_m}$ as in \eqref{eq:unshifted-trunc}, $Q_j \in \mathbb{R}^{(m+1)n \times j}$ and $\underline{I}_j$ as in \eqref{eq:Im} and $\underline{H}_j \in \mathbb{R}^{(j+1) \times j}$. The matrix $M_m K_m^{-1}$ is given by
\begin{align} \label{eq:MK-inv}
	M_m K_m^{-1} =
	\begin{bmatrix}
	0 & 0 & 0 & \cdots & 0 \\
	A_{0}^{-1}& -A_{0}^{-1}A_1 & -A_{0}^{-1} A_2 & \cdots & -A_{0}^{-1}A_m \\
	& I & & & \\
	& & \ddots & & \\
	& & & I & 0
	\end{bmatrix} \in \mathbb{R}^{(m+1)n \times (m+1)n}.
\end{align}
Consider first the action of the truncated matrix product $M_m K_m^{-1}$. From \eqref{eq:multMK} it is clear that in forming the Arnoldi factorization we must apply $A_0^{-1}$ to a new vector of size $n$ on each iteration.

In a flexible GMRES setting, we apply an approximation of the preconditioning matrix on each iteration, denoted $\tilde{K}_{m_i}^{-1}$, $i=1,\ldots,j.$ We refer to this particular setting as inexact GMRES. The inexactness here refers to approximating the most costly component of the preconditioner, i.e., the action of $A_0^{-1}$ on a different vector each time it is applied. This approximation, denoted $\tilde{A}_{0}^{-1}$, is permitted to change from one iteration to the next.\footnote{$\tilde{A}_0^{-1}$ is implicitly defined by how we approximate the action of $A_0^{-1}$ on each iteration.} The modification described here does not change the nonzero structure of the matrix in \eqref{eq:MK-inv}. The shifted inexact Arnoldi factorization on iteration $j$ is given by
\begin{align} \label{eq:shift-in-trunc-arn}
	Q_j - \mu M_m \tilde{Z}_j = Q_{j+1} (\underline{I}_j - \mu \underline{H}_j),
\end{align}
where
\begin{align} \label{eq:flex-mult}
	\tilde{z}_i = \tilde{K}^{-1}_{m_i} q_i,
\end{align}
$\tilde{Z}_j = \begin{bmatrix} \tilde{z}_1, \ldots, \tilde{z}_j \end{bmatrix} \in \mathbb{R}^{(m+1)n \times j}$, and $M_m$ as in \eqref{eq:unshifted-trunc}. The following theorem characterizes the important structure of the matrices $Q_i$ and $\tilde{Z}_i$, which will be proven as a consequence of the more technical result Lemma~\ref{matvec}.
\begin{theorem} \label{lemma1}
If we apply inexact GMRES to the linearized system in \eqref{eq:unshifted-trunc} with preconditioner $\tilde{K}_{m_i}^{-1}$, the nonzero structure in both $Q_i$ and $\tilde{Z}_i$ as in \eqref{eq:shift-in-trunc-arn} grows by one block in each iteration $i=1,\ldots,j$.
\end{theorem}

Note, the result of Theorem~\ref{lemma1} holds if $A_0^{-1}$ is applied exactly. We show in the following Lemma that the action of the operator $M_m \tilde{K}_{m_i}^{-1}$ (or $M_m K_m^{-1}$) can be computed efficiently. Proofs of Lemma~\ref{matvec} and Theorem~\ref{lemma1} respectively follow after.
\begin{lemma} \label{matvec}
Let $M_m \tilde{K}_{m_i}^{-1}$ be as in \eqref{eq:MK-inv} with approximation $\tilde{A}_0^{-1}$ on iteration $i$. For $W = \begin{bmatrix} w_0, w_1, \ldots, w_m \end{bmatrix} \in \mathbb{R}^{n \times (m+1)}$,
\[
    M_m \tilde{K}_{m_i}^{-1} \vectorize(W) = \vectorize(\bar{0}, \tilde{w}, W \tilde{D}_m^T) \in \mathbb{R}^{n(m+1)},
\]
where
\begin{equation} \label{eq:tilde-w-small}
    \tilde{w} = \tilde{A}_0^{-1} w_0 -\tilde{A}_0^{-1} \left(\sum_{\ell=1}^{m} A_{\ell} w_{\ell} \right) \in \mathbb{R}^{n},
\end{equation}
$\bar{0} \in \RR^{n}$, and $\tilde{D}_m$ is as in \eqref{eq:tilde-D-m}.
\end{lemma}
{\em Proof.}
Lemma~\ref{matvec} is proved as follows. We rewrite the product between the matrix in \eqref{eq:MK-inv} with approximation $\tilde{A}_0^{-1}$ and $\vectorize(W)$ as
\[
	M_m \tilde{K}_{m_i}^{-1} \vectorize(W) =
	\begin{bmatrix}
	\left(e^T \otimes D_0 \right) \vectorize(W)  \\
	\left(\tilde{A}_0^{-1} \left[ I, -A_1, -A_2, \ldots, -A_m \right]\right) \vectorize(W) \\
	\left(\tilde{D}_m \otimes I\right)\vectorize(W)
	\end{bmatrix},
\]
where $e=\begin{bmatrix}1,\ldots,1\end{bmatrix}^T \in \RR^{m+1}$, $D_0 \in \RR^{n \times n}$ is the zero matrix, and $\tilde{D}_m$ is as in \eqref{eq:tilde-D-m}. The first block row reduces to
\[
	\vectorize(D_0 W e) = \begin{bmatrix} 0, \ldots, 0\end{bmatrix}^T \in \RR^{n}
\]
and the second block row simplifies as
\[
	\tilde{A}_0^{-1} \left( w_0 -A_1 w_1 -A_2 w_2 -\ldots -A_m w_m \right) = \tilde{A}_0^{-1} w_0 - \tilde{A}_0^{-1} \sum_{\ell=1}^m A_{\ell} w_{\ell} \in \RR^{n}.
\]
The remaining product is computed as
\[
	\left(\tilde{D}_m \otimes I\right)\vectorize(W) = \vectorize(W \tilde{D}_m^T) = \vectorize( \begin{bmatrix} w_1,\ldots,w_{m-1} \end{bmatrix}) \in \RR^{n(m-1)},
\]
where we have used that $(X^T \otimes Y) \vectorize(Z)= \vectorize(YZX)$. \hfill $\endproof$

Note that if $W$ is of the form
\begin{align} 
	W = \begin{bmatrix} \hat{W},0, \dots ,0 \end{bmatrix} \in \mathbb{R}^{n \times (m+1)}, \text{ }
	\hat{W} = \begin{bmatrix} w_0, \ldots, w_k \end{bmatrix} \in \mathbb{R}^{n \times (k+1)}
\end{align}
with $k < m$, then by a direct consequence of Lemma~\ref{matvec}, the matrix-vector product reduces to
\begin{align} \label{eq:smaller-prod}
	M_m \tilde{K}_{m_i}^{-1} \vectorize{(W)} = \vectorize{(\bar{0}, \tilde{w}, \hat{W} \tilde{D}_{k}^T,0, \dots ,0)} \in \mathbb{R}^{n(m+1)},
\end{align}
where $\tilde{w}$ is as in \eqref{eq:tilde-w-small} with $m$ replaced by $k$.

{\em Proof.}
Theorem~\ref{lemma1} is proved by the following reasoning. Take $q_1=c_m/\norm{c_m}$ with $c_m$ as in \eqref{eq:unshifted-trunc}. A direct result of Lemma~\ref{matvec} is the structure of the resulting vector is changed only by an expanding the number of nonzero entries in the first block, and the tail of zeros is preserved after orthogonalization (structure of $Q_i$). The same result holds for $\tilde{Z}_i$, where we require the matrix-vector product in \eqref{eq:flex-mult} on each iteration. The preconditioning matrix has a block diagonal structure and thus the nonzero structure is preserved after multiplication. \hfill $\endproof$

To visualize Theorem~\ref{lemma1} and Lemma~\ref{matvec}, consider the matrix-vector products $M_m \tilde{K}_{m_i} q_i$ in the Arnoldi factorization in \eqref{eq:shift-in-trunc-arn}, started on the vector $q_1 = c_m/ \norm{c_m}$ as in \eqref{eq:unshifted-trunc}. The first product is given by
\begin{align} \label{eq:first-product} \footnotesize
	\begin{bmatrix}
	0 & 0 & 0 & \cdots & 0 \\
	\tilde{A}_{0}^{-1}& -\tilde{A}_{0}^{-1}A_1 & -\tilde{A}_{0}^{-1} A_2 & \cdots & -\tilde{A}_{0}^{-1}A_m \\
	& I & & & \\
	& & \ddots & & \\
	& & & I & 0
	\end{bmatrix}
	\begin{bmatrix}
	b/\norm{b} \\
	0 \\
	0 \\
	\vdots \\
	0
	\end{bmatrix} =
	\begin{bmatrix}
	\bar{0} \\
	\tilde{A}_0^{-1} b/\norm{b} \\
	0 \\
	\vdots \\
	0
	\end{bmatrix},
\end{align} \normalsize
and subsequent multiplications by
\begin{align} \label{eq:star} \footnotesize
	\begin{bmatrix}
	0 & 0 & 0 & \cdots & \cdots & \cdots & \cdots & 0 \\
	\tilde{A}_{0}^{-1}& -\tilde{A}_{0}^{-1}A_1 & -\tilde{A}_{0}^{-1} A_2 & \cdots & \cdots & \cdots & \cdots & -\tilde{A}_{0}^{-1}A_m \\
	& I & & & & & & \\
	& & \ddots & & & & & \\
	& & & & & & & \\
	& & & &  & & & \\
	& & & & & \ddots & & \\
	& & & & & & I & 0
	\end{bmatrix}
	\begin{bmatrix}
	\bar{0} \\
	* \\
	\vdots \\
	* \\
	0 \\
	0 \\
	\vdots \\
	0
	\end{bmatrix} =
	\begin{bmatrix}
	\bar{0} \\
	* \\
	\vdots \\
	* \\
	* \\
	0 \\
	\vdots \\
	0
	\end{bmatrix},
\end{align} \normalsize
where $\bar{0} \in \mathbb{R}^n$ and we have used that the tail of zeros is preserved after orthogonalization. Here the entries marked with $*$ represent nonzero blocks of dimension $n$. Note the number of nonzero blocks grows by one on each iteration of the Arnoldi factorization. More specifically, on the $i$th iteration we have nonzero blocks in positions $2$ through $i$ in the solution vector in equations \eqref{eq:first-product} and \eqref{eq:star}.

The above process is summarized in Algorithm~\ref{alg:FIGMRES}, where we return a parametrized solution $\tilde{x}_j(\mu)$ as an approximation to \eqref{eq:our-prob}. We evaluate this solution for many values of $\mu$ without executing the algorithm again. Note that within Algorithm~\ref{alg:FIGMRES}, $q_i$ refers only to the nonzero elements above the tail of zeros. The nonzero elements of $\tilde{z}_i$, considering an inexact application of the preconditioner as in \eqref{eq:flex-mult}, can be expressed as
 \begin{align} \label{eq:z-tilde-for-alg}
 	\vectorize{(\tilde{w}, \hat{W} \tilde{D}_i)} \in \mathbb{R}^{in},
 \end{align}
where $\vectorize{(\hat{W})} = \begin{bmatrix} q_i^T, \bar{0}^T \end{bmatrix}^T \in \mathbb{R}^{(i+1)n}$ with $q_i \in \mathbb{R}^{in}$, $\bar{0} \in \mathbb{R}^n$, $\tilde{w}$ as in \eqref{eq:tilde-w-small} and $\tilde{D}_i$ as in \eqref{eq:tilde-D-m}. The elements of $y = M \tilde{z}_i$ above the tail of zeros are given by
 \begin{align} \label{eq:y-for-alg}
 	\begin{bmatrix} \bar{0}^T, \tilde{z}_i^T \end{bmatrix}^T \in \mathbb{R}^{(i+1)n}
 \end{align}
and the orthogonalization on line \ref{alg:gs-line1} is performed with a (possibly repeated) Gram--Schmidt procedure. Changing the last line of Algorithm~\ref{alg:FIGMRES} to
\begin{align} \label{eq:change_to_v}
    \tilde{v}_j (\mu): \mu \mapsto \tilde{Z}_j \argmin_{w \in \mathbb{R}^j} \norm{ \left( \underline{I}_j - \mu \underline{H}_j \right) w - e_1 \norm{c_m} }_2,
\end{align}
where $\tilde{v}_j(\mu) \in \mathbb{R}^{(j+1)n}$ results in an approximation to the truncated system in \eqref{eq:unshifted-trunc}.
\begin{remark}[Storing nonzero elements]
On iteration $i$ of Algorithm~\ref{alg:FIGMRES}, only the first $n(i+1)$ elements in the resulting matrix-vector products are stored. This is in contrast to the standard inexact GMRES algorithm, where we require the storage of full vectors.
\end{remark}
\begin{remark}[Inexact application of the preconditioner]
In Section~\ref{sec:inner-tol} we quantify how inexactly the preconditioner can be applied, without disrupting the convergence of Algorithm~\ref{alg:FIGMRES}. The specific criteria indicated on Line~\ref{alg:Kinv} of the algorithm, along with input variables $\varepsilon$ and $\ell_j$, are explained in further detail in this section.
\end{remark}

\begin{algorithm}[h!]
\SetAlgoLined
\SetKwInOut{Input}{Input}\SetKwInOut{Output}{Output}
\SetKwFor{FOR}{for}{do}{end}
\Input{$j \in \mathbb{Z}^{+}$, desired dimension of the Krylov subspace, \\ $A_{\ell}$, $\ell = 0,1,\ldots,j$ as in \eqref{eq:taylor-coeff}, \\ $b \in \mathbb{R}^n$, \\ $\mu \in \mathbb{R}$, $\ell_j \in \RR$, $\varepsilon > 0$}
\Output{Approximate solution $\tilde{x}_j(\mu)$ to \eqref{eq:our-prob}}
$\tilde{Z}_0$ = empty matrix \\
$Q_1 = b/ \norm{b}$ \\
\FOR{$i = 1,2, \dots, j$}{
Compute nonzero blocks of $\tilde{z}_i$ as in \eqref{eq:z-tilde-for-alg} s.t. \eqref{eq:pj-ineq} holds \label{alg:Kinv} \\
Compute elements of $y$ above tail of zeros as in \eqref{eq:y-for-alg} \\
Orthogonalize $y$ against $q_1,\ldots,q_i$ by a Gram--Schmidt process \label{alg:gs-line1}  \\
    \nonl \Indp $h_i = Q_i^{T}y$ \\
    \nonl $y_{\perp} = y - Q_i h_i$ \\
    \Indm
    	Possibly repeat Step~\ref{alg:gs-line1} \\
    	Compute $\beta_i = ||y_{\perp}||$ \\
	Let $q_{i+1} = y_{\perp}/\beta_i$ \\
	Expand $Q_i$ by $n$ rows of zeros \\
	Append $Q_i$ s.t. $Q_{i+1} = [Q_i, q_{i+1}]$ \\
	Append $\tilde{Z}_{i-1}$ s.t. $\tilde{Z}_i = [\tilde{Z}_{i-1}, \tilde{z}_i]$ \\
	Let
	$\underline{H}_i =
    	\begin{bmatrix}
    	\underline{H}_{i-1} & h_i \\
    	0 & \beta_i
    	\end{bmatrix} \in \mathbb{R}^{(i+1) \times i}$
}
Return function handle
\begin{align*} 
    \tilde{x}_j (\mu): \mu \mapsto \tilde{Z}_j (1:n,:)  \argmin_{w \in \mathbb{R}^j} \norm{ \left( \underline{I}_j - \mu \underline{H}_j \right) w - e_1 \norm{b} }_2
\end{align*}
with $\underline{I}_j$ as in \eqref{eq:Im}.
\caption{Inexact infinite GMRES}
\label{alg:FIGMRES}
\end{algorithm}

\subsubsection{Extension to infinity} \label{sec:ext-to-infty}
Note that Algorithm~\ref{alg:FIGMRES} is independent of the truncation parameter $m$, shown in Lemma~\ref{matvec}. More precisely, we proved the number of operations to calculate the product in \eqref{eq:smaller-prod} is independent of the length of the vector, assuming a starting vector with a fixed number of nonzero blocks. Hence, in theory we can take $m \rightarrow \infty$ and compute the product of a vector with an infinite tail of zeros with an infinite companion matrix in a finite number of linear algebra operations. As in the proof of Theorem~\ref{lemma1}, the infinite tail of zeros on each of the resulting vectors are preserved even after orthogonalization. Hence, Algorithm~\ref{alg:FIGMRES} is equivalent to the application of standard inexact GMRES for parameterized linear systems applied to the infinite matrices, i.e., corresponding to $m \rightarrow \infty$. This procedure is completely analogous to the infinite Arnoldi method \cite{Jarlebring:2012:INFARNOLDI}.

The process described above therefore represents the Taylor series of $A(\mu)$ as in \eqref{eq:our-prob} without any truncation error, i.e., Algorithm~\ref{alg:FIGMRES} as written indeed approximates the solution in \eqref{eq:lin} with $m \rightarrow \infty$ and equivalently, the solution in \eqref{eq:our-prob}.

\subsection{Dynamic inner stopping criteria} \label{sec:inner-tol}
In Algorithm \ref{alg:FIGMRES}, we have adapted the inexact preconditioning strategy to the infinite GMRES method applied to \eqref{eq:our-prob};
however, we have not specified a stopping strategy for the inner iteration.  We could simply specify a fixed, appropriately-stringent tolerance to obtain
convergence with no degredation due to the inexact preconditioning.  Instead, in this section, we adapt the theory of inexact Krylov subspace
methods \cite{inexKry} with dynamic stopping criteria to show that we are able to relax the inner convergence tolerance such that
we do not degrade convergence of the outer iteration. This framework allows for a more efficient algorithm without sacrificing accuracy. We adapt this criteria to the application of the inexact preconditioner in Algorithm~\ref{alg:FIGMRES}.

In \cite{inexKry}, the authors study the convergence properties of a GMRES iteration in which the matrix-vector product (preconditioned or
unpreconditioned) is applied inexactly.  This is of great interest since many matrix-free operators are applied implicitely as procedures that are
often computationally expensive to apply with high accuracy and are indeed able to be tuned for application at lower accuracy.  The question the
authors sought to answer was: how accurately must one apply the inexact matrix-vector product at a given iteration such that the convergence of
GMRES is not degraded?  It is shown that it is sufficient to choose the inner error tolerance to be roughly inversely proportional to the size
of the GMRES residual. This means that as the method converges, the inexact matrix-vector product tolerance can be relaxed and thus becomes
cheaper as iteration proceeds.

In this work, we are employing an inexact preconditioning strategy; so it is the preconditioned matrix-vector product which is inexact, and we adapt
the results of in \cite{inexKry} to determine a criteria for how accurately the inexact precondtioner must be applied.
More specifically, at the $j$th iteration of Algorithm~\ref{alg:FIGMRES} we form the inexact Arnoldi factorization \eqref{eq:shift-in-trunc-arn}. This is achieved by computing $\tilde{z}_i$ for $i=1,\ldots,j$ as in \eqref{eq:flex-mult} via an inexact application of $A_0^{-1}$. The application of $A_0^{-1}$ by an iterative method can be performed to different tolerance level at each Arnoldi step. The following analysis shows how the accuracy of the application of $A_0^{-1}$ can be chosen dynamically.

Consider solving \eqref{eq:our-prob} for a fixed $\mu$ using Algorithm \ref{alg:FIGMRES}.
Let $\tilde{v}_j (\mu) \in \mathbb{R}^{(j+1) \times n}$ be as in \eqref{eq:change_to_v}, i.e., the $j$th parametrized approximation returned from Algorithm~\ref{alg:FIGMRES} applied to the system in \eqref{eq:unshifted-trunc}.  Assume $m \to \infty$, i.e., the system in \eqref{eq:unshifted-trunc} contains an infinite Taylor series expansion without truncation (see Section~\ref{sec:ext-to-infty}). At iteration $j$, the GMRES minimization becomes
\begin{subequations} \label{eq:two-eq}
\begin{eqnarray}
	w_j &=& \argmin_{w \in \mathbb{R}^j} \norm{(\underline{I}_j - \mu \underline{H}_j)w - e_1 \norm{c_m}}_2 \label{eq:wm}, \\
	\tilde{v}_j(\mu) &=& \tilde{Z}_j w_j
\end{eqnarray}
\end{subequations}
with $\tilde{Z}_j, \underline{H}_j$ as in \eqref{eq:shift-in-trunc-arn}. The inner residual on the $i$th Arnoldi iteration, denoted $p_i$, represents how inexactly we apply the preconditioning matrix $K_m^{-1}$, i.e.,
\begin{align} \label{eq:pj1}
	p_i = K_m \tilde{z}_i - q_i,
\end{align}
and thus in matrix form, we have
\begin{align} \label{eq:res-mat}
	K_m \tilde{Z}_j = Q_j + P_j,
\end{align}
where $P_j = \begin{bmatrix} p_1, \ldots, p_j \end{bmatrix}$ after $j$ iterations.

Thus, using \eqref{eq:shift-in-trunc-arn} and \eqref{eq:res-mat}, we express $r_j$, the residual of the $j$th approximate solution to \eqref{eq:unshifted-trunc} as follows.
\begin{subequations}\label{eq:r-rtilde}
\begin{eqnarray}
	r_j &=&  c_m - (K_m - \mu M_m) \tilde{v}_j (\mu) \\
	&=& c_m - (K_m - \mu M_m) \tilde{Z}_j w_j \label{eq:line2} \\
	&=& c_m - K_m \tilde{Z}_j w_j + \mu M_m \tilde{Z}_j w_j \\
	&=& c_m - K_m \tilde{Z}_j w_j - \big( Q_{j+1} (\underline{I}_j - \mu \underline{H}_j) - Q_j \big) w_j \\
	&=& c_m - (Q_j + P_j ) w_j - Q_{j+1} (\underline{I}_j - \mu \underline{H}_j) w_j + Q_j w_j \\
	&=& Q_{j+1} \big( e_1 \norm{c_m} - (\underline{I}_j - \mu \underline{H}_j) w_j \big) - P_j w_j \\
	&=& \tilde{r}_j - P_j w_j \label{eq:separate}
\end{eqnarray}
\end{subequations}
We refer to the vector $\tilde{r}_j$ above as the \emph{exact residual}. This vector expresses the residual at the $j$th iteration under the assumption that the preconditioning matrix $K_m^{-1}$ is applied exactly on each iteration. Equation \ref{eq:separate} enables us to quantify the
error introduced from the inexact preconditioning from the residual of the GMRES applied to the shifted system \eqref{eq:shifted-trunc},
which we define as
\begin{align} \label{eq:distance}
	\delta_j := \norm{r_j - \tilde{r}_j} = \norm{P_j w_j}.
\end{align}
If $\delta_j$ is small, then inexact GMRES has produced
an approximation close to that produced by GMRES with exact preconditioning.

The goal is to quantify \emph{how} inexactly we can apply
the preconditioner at each iteration $j$ and guarantee that $\delta_j$ is sufficiently small.
The following proposition quantifies this by specifying how large the norm of \eqref{eq:pj1} can be
such that we can guarantee a sufficiently small value for the quantity $\delta_j$.

\begin{proposition} \label{prop1}
Let $\varepsilon > 0$ and $w_j = \begin{bmatrix} \eta_1^{(j)}, \ldots, \eta_j^{(j)} \end{bmatrix}^T \in \mathbb{R}^j$ as in \eqref{eq:wm}. Let the residual of Algorithm~\ref{alg:FIGMRES} applied to \eqref{eq:unshifted-trunc} on iteration $j$ be given as in \eqref{eq:line2} and the exact residual
\[
	\tilde{r}_i = Q_{i+1} \big( e_1 \norm{c_m} - (\underline{I}_i - \mu \underline{H}_i) w_i \big)
\]
for $i=1,\ldots,j$. If at each iteration $i \leq j$ the inner residual $p_i$ defined in \eqref{eq:pj1}
satisfies
\begin{align} \label{eq:stop-crit}
	\norm{p_i} \leq \frac{ \sigma_j(\underline{I}_j - \mu \underline{H}_j)}{j} \frac{1}{\norm{\tilde{r}_{i-1}}} \varepsilon =: \varepsilon_{\text{inner}}^{(i)},
 \end{align}
then $\delta_j \leq \varepsilon$, where $\delta_j$ as in \eqref{eq:distance}.
\end{proposition}

{\em Proof.}
A direct application of \cite[Lemma~5.1 and Proposition 9.1]{inexKry} gives the desired result. \hfill $\endproof$

We note that we can express \eqref{eq:stop-crit} as
\begin{align} \label{eq:pj-ineq}
	\norm{p_i} \leq \ell_j \frac{1}{\norm{\tilde{r}_{i-1}}} \varepsilon = \varepsilon_{\text{inner}}^{(i)}
\end{align}
for $i=1,\ldots,j$, where $\ell_j$ is a constant independent of $i$. Assuming $\varepsilon > 0$ is fixed and sufficiently small, this
indicates that the quantities $\norm{p_i}$ and $\norm{\tilde{r}_{i-1}}$ are inversely related since
\begin{align} \label{eq:inverse}
	\norm{p_i} \norm{\tilde{r}_{i-1}} \approx \ell_j \varepsilon.
\end{align}
Thus, \eqref{eq:pj-ineq} admits the interpretation: the smaller the norm of the exact residual becomes, the more inexactly we can apply preconditioner $K_m^{-1}$ as in \eqref{eq:flex-mult}. In practice, this corresponds to applying the matrix $A_0^{-1}$ via an iterative method with an increasingly relaxed
tolerance as the residual norm decreases.

\section{Convergence theory} \label{sec:convergence-the}
Our main interest is to apply Algorithm~\ref{alg:FIGMRES} to approximate \eqref{eq:our-prob} for many different values of $\mu$. We therefore examine the relationship between the performance of the algorithm and the magnitude of the parameter $\mu$. Consider the norm of the residual where we apply Algorithm~\ref{alg:FIGMRES} to \eqref{eq:unshifted-trunc}, i.e., at the $j$th iteration we compute
\begin{align} \label{eq:big-res}
	\norm{r_j} = \norm{c_m - (K_m-\mu M_m) \tilde{v}_j(\mu)}
\end{align}
with $K_m$, $M_m$ and $c_m$ as in \eqref{eq:unshifted-trunc} and $\tilde{v}_j(\mu)$ is the $j$th iterate of Algorithm~\ref{alg:FIGMRES} as in \eqref{eq:change_to_v}. We can express $\norm{r_j}$ in terms of the exact residual and an error term as in \eqref{eq:r-rtilde}, i.e.,
\begin{subequations}
\begin{eqnarray*}
	\norm{r_j} &=& \norm{Q_{j+1} \big( e_1 \norm{c_m} - (\underline{I}_j - \mu \underline{H}_j) w_j \big) - P_j w_j}, \\
	&=& \norm{\tilde{r}_j - P_j w_j},
\end{eqnarray*}
\end{subequations}
where the matrices $Q_j$, $H_j$ are as in \eqref{eq:shift-in-trunc-arn}, and the matrix $P_j$ contains the inner residual vectors as in \eqref{eq:pj1}. Thus,
\begin{align} \label{eq:r-ineq}
	\norm{r_j} \leq \norm{\tilde{r}_j} + \norm{P_j w_j}.
\end{align}
The following theorem gives a bound on $\norm{r_j}$ based on this inequality.

\begin{theorem} \label{eq:theorem-res}
Suppose Algorithm~\ref{alg:FIGMRES} is executed with an inexact solve satisfying the criteria \eqref{eq:stop-crit}.
Then, the inexact residual vector at iteration $j$ satisfies
\begin{align} \label{eq:r-eq}
\|r_j\|\le \min_{p \in \mathcal{P}_j^0} \norm{p (I - \mu M_m K_m^{-1}) c_m} + \varepsilon,
\end{align}
where
\[
	\mathcal{P}_j^0 := \{ \text{polynomials } p \text{ of degree} \leq j \text{ with } p(0) = 1 \}.
\]
\end{theorem}
{\em Proof.}
We consider first the norm of the exact residual $\norm{\tilde{r}_i}$ for each $i=1,2,\ldots, j$.
Each exact residual is associated to an \emph{exact approximation}, i.e., there exists $\tilde{u}_i(\mu)$ such that
\begin{align} \label{eq:u-tilde-i}
	\tilde{r}_i = c_m - (K_m - \mu M_m) \tilde{u}_i(\mu).
\end{align}
		We express $\norm{\tilde{r}_j}$ as
	\begin{subequations} \label{eq:calc-r-tilde}
	\begin{eqnarray}
		\norm{\tilde{r}_j} &=& \norm{c_m - (K_m - \mu M_m) \tilde{u}_j (\mu)} \label{eq:norm-r-tilde} \\
		&=& \norm{c_m - (K_m - \mu M_m) K_m^{-1} K_m \tilde{u}_j (\mu)} \\
		&=& \norm{c_m - (I-\mu M_m K_m^{-1}) \hat{\tilde{u}}_j (\mu)},
	\end{eqnarray}
	\end{subequations}
where $\hat{\tilde{u}}_j(\mu) = K_m \tilde{u}_j (\mu)$.
The rest of the proof follows
from applying the standard convergence theory of GMRES \cite{Saad:1986:GMRES} to the
first term in \eqref{eq:u-tilde-i}, combined with applying Proposition~\ref{prop1} to the second term. \hfill $\endproof$

%
%

To illustrate a practical implication of Theorem~\ref{eq:theorem-res}, we further bound \eqref{eq:r-eq} by choosing a particular polynomial,
\begin{subequations}
\begin{eqnarray}
	\norm{r_j} &\leq& \min_{p \in \mathcal{P}_j^0} \norm{p (I - \mu M_m K_m^{-1})c_m} + \varepsilon \\
	&\leq& \norm{q (I - \mu M_m K_m^{-1})c_m} + \varepsilon, \label{eq:r-tilde-2}
\end{eqnarray}
\end{subequations}
where the polynomial $q$ is defined as
\[
	q(z) := \frac{(1-z)^{j-k} \prod_{i=1}^k (1- z - \mu \gamma_i)}{\prod_{i=1}^k (1-\mu \gamma_i)},
\]
and $\gamma_{k+1}$ is the ($k+1$)th largest eigenvalue in modulus of $M_m K_m^{-1}$. Here $q$ is a normalized degree $j$ Zarantonello polynomial \cite{ZaranVarga}. The inequality in \eqref{eq:r-tilde-2} can thus be rewritten as
\begin{align} \label{eq:r-tilde-bound}
	\norm{r_j} \leq \norm{\frac{(\mu M_m K_m^{-1})^{j-k} \prod_{i=1}^k (\mu M_m K_m^{-1} - \mu \gamma_i)}{\prod_{i=1}^k (1-\mu \gamma_i)}} \norm{c_m} + \varepsilon,
\end{align}
and we obtain
\begin{align} \label{eq:conv-bound}
	\norm{r_j} \leq ( |\mu| |\gamma_{k+1}| )^j + \varepsilon
\end{align}
for sufficiently large $j$ (see \cite[Theorem~4.1]{JarlebringCorrenty1}). This indicates that for a sufficiently small neighborhood $\mathcal{N}\subset\mathbb{R}$ around the origin,
we expect that for $\mu_1,\mu_2\in\mathcal{N}$ with $|\mu_1|<|\mu_2|$ that the convergence will be more rapid for the residual associated to
$\mu_1$.

As noted in Section~\ref{sec:ext-to-infty}, as $m \rightarrow \infty$, Algorithm~\ref{alg:FIGMRES} as stated approximates the solution to \eqref{eq:lin} without truncation, i.e., it approximates the solution to \eqref{eq:our-prob} without neglecting any of the terms in the infinite Taylor series expansion \eqref{eq:A-taylor}. Therefore, the approximation error comes entirely from the iterative method. In practice, Algorithm~\ref{alg:FIGMRES} converges fast for small values of $\mu$.

\begin{remark}[Different $\mu$ in evaluation and stopping criteria]
Theorem 4.1 enables us to characterize the convergence of Algorithm~1 in terms of
$\mu$ and the	 spectra of $M_{m}K_{m}^{-1}$.
Let $\mathcal{N}$ be the aforementioned neighborhood.
Consider solving the problem using Algorithm \ref{alg:FIGMRES} for many different  parameters
$\mu\in\left\lbrace \mu_1,\mu_2,\ldots,\mu_k \right\rbrace\subset\mathcal{N}$.
We need to use the residual associated to a specific $\mu$ for the stopping criteria.
This means we must choose a specific $\mu$, and the residual for this parameter is the one to which we apply the convergence theory.
Suppose that we index the parameters such that
\begin{align*}
	|\mu_1|\leq|\mu_2|\leq\cdots\leq|\mu_k|.
\end{align*}
Since all parameters are in $\mathcal{N}$,
the above theory suggests that Algorithm \ref{alg:FIGMRES} converges more quickly for smaller values of $\mu\in\mathcal{N}$, which then means
that the bound \eqref{eq:pj-ineq} becomes less strict more quickly.  This implies that the dynamic inner tolerance for the inexact
preconditioner should be chosen using the residual associated to $\mu_k$, as this will guarantee that the dynamic inner tolerances will
also work for all other parameters in the list.
\end{remark}

\section{Simulations} \label{sec:simulations}
\subsection{Illustration of the theory} \label{sec:time-delay}
The simulations in this work were carried out on a system with a 2.3 GHz Dual-Core Intel Core i5 processor and 16 GB RAM. The software for all examples in this section are available online.\footnote{\url{https://github.com/siobhanie/FlexInfGMRES}}

Consider approximating the solution to \eqref{eq:our-prob} with Algorithm~\ref{alg:FIGMRES}, where
\begin{align} \label{eq:delay-A}
	A(\mu) = -\mu I + A_0 + A_1 e^{-\mu}
\end{align}
with $A_0, A_1 \in \RR^{n \times n}$ random banded matrices and $b \in \RR^{n}$ random. The solution to this system is the transfer function of the time-delay system described by
\begin{subequations}
\begin{eqnarray*}
	\dot{x}(t) &=& A_0 x(t) + A_1 x (t - \tau) - b u(t), \\
	y(t) &=& C^T x(t),
\end{eqnarray*}
\end{subequations}
where entire state is the output, i.e. $C=I \in \mathbb{R}^{n \times n}$. Specifically, the transfer function is obtained by applying the Laplace transform to the state equation under the condition $x(0)=0$ and $\mu$ is the Laplace variable (see standard references for time-delay
systems \cite{Gu:2003:STABILITY,Michiels2011,Michiels:2007:STABILITYBOOK}). Here the vector $b \in \mathbb{R}^n$ is the external force, $x(t) \in \mathbb{R}^n$ is the state vector, $u(t)$ is the input, $y(t)$ is the output, and $\tau > 0$ is the delay. Assume without loss of generality that $\tau = 1$.

In the following simulation, we replace line \ref{alg:Kinv} of Algorithm~\ref{alg:FIGMRES} with
\begin{align} \label{eq:z-pert}
	\tilde{z}_i = K_m^{-1} q_i + \Delta_i,
\end{align}
where $\Delta_i$ is random vector. The above introduces error to line \ref{alg:Kinv} of Algorithm~\ref{alg:FIGMRES} in a slightly different way than considered previously. Here the preconditioning matrix $K_m^{-1}$ as in \eqref{eq:Km-inv} is applied exactly but we introduce error with a random vector. The magnitude of the random vector $\Delta_i$ is determined in a way analogous to Proposition~\ref{prop1}. More precisely, we consider a bound on $p_i$ as in \eqref{eq:pj-ineq} where $p_i$ is the inner residual vector as in \eqref{eq:pj1} and $\ell_j = 1$. In the simulations which follow, we visualize the relative residual
\begin{align} \label{eq:rel-res-rand}
	\frac{\norm{A(\mu)\tilde{x}_i(\mu) - b}}{\norm{b}},
\end{align}
where $\tilde{x}_i$ the output of Algorithm~\ref{alg:FIGMRES} on iteration $i$.

\begin{figure}[h]
	\centering
	\includegraphics{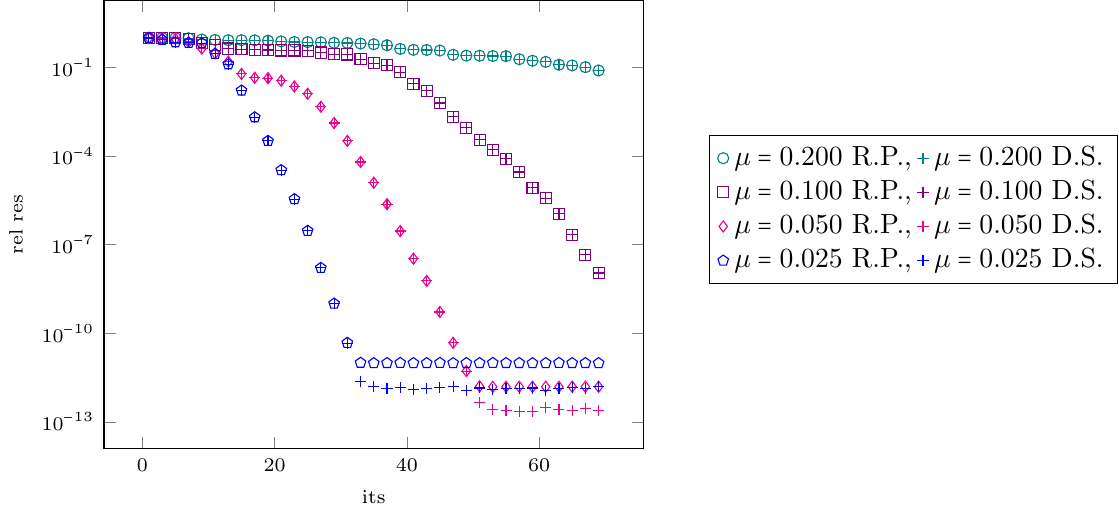}

	\caption{Convergence of Algorithm~\ref{alg:FIGMRES} for solving \eqref{eq:unshifted-trunc}, $m \rightarrow \infty$. Here $A(\mu)$ as in \eqref{eq:delay-A} with $n=1000$, random perturbations (R.P.) as in \eqref{eq:z-pert}, relative residual as in \eqref{eq:rel-res-rand}, $\ell_j = 1$, $\varepsilon = 10^{-10}$, and $\varepsilon_{\text{inner}}^{(i)}$ as in \eqref{eq:pj-ineq} correspond to largest $\mu$ (not plotted). We compare to direct solve (D.S.) without random perturbations on each iteration.}
	\label{fig:asymfigs}
\end{figure}

\begin{figure}
	\begin{subfigure}[T]{0.3\columnwidth}
	\includegraphics{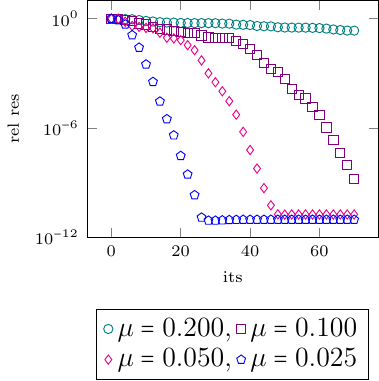}

	\caption{$\varepsilon=10^{-10}$}
	\label{fig:rtilde1}
	\end{subfigure}
	\hfill
	\begin{subfigure}[T]{0.3\columnwidth}
	
	\includegraphics{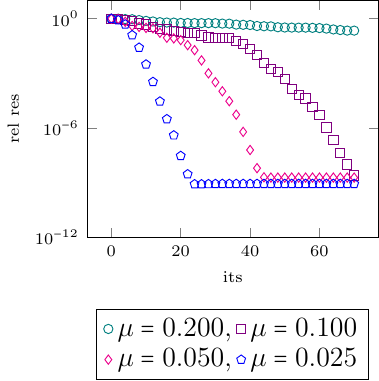}

	\caption{$\varepsilon=10^{-8}$}
	\label{fig:rtilde2}
	\end{subfigure}
	\hfill
	\begin{subfigure}[T]{0.3\columnwidth}
	
	\includegraphics{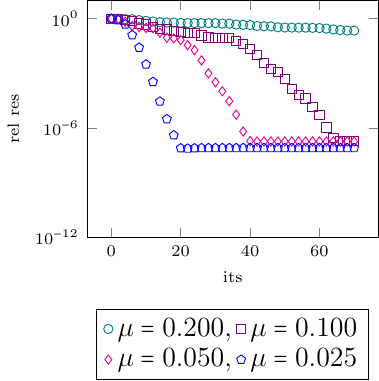}

	\caption{$\varepsilon=10^{-6}$}
	\label{fig:rtilde3}
	\end{subfigure}
	\caption{Convergence of Algorithm~\ref{alg:FIGMRES} for solving \eqref{eq:unshifted-trunc}, $m \rightarrow \infty$. Here $A(\mu)$ as in \eqref{eq:delay-A} with $n=1000$, random perturbations as in \eqref{eq:z-pert}, norm of relative residual as in \eqref{eq:rel-res-rand}, $\ell_j = 1$ for different values of $\varepsilon$ as in \eqref{eq:pj-ineq}.}
	\label{fig:rtilde}
\end{figure}

\begin{figure}
	\begin{subfigure}[T]{0.3\columnwidth}
	
	\includegraphics{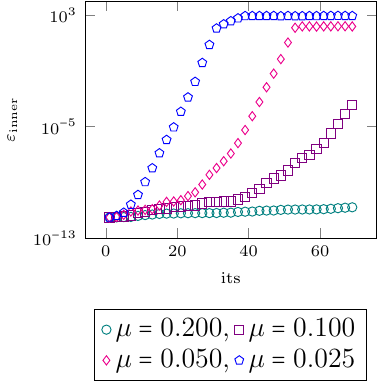}

		\caption{$\varepsilon=10^{-10}$}
		\label{fig:einner1}
	\end{subfigure}
	\hfill
	\begin{subfigure}[T]{0.3\columnwidth}
	
	\includegraphics{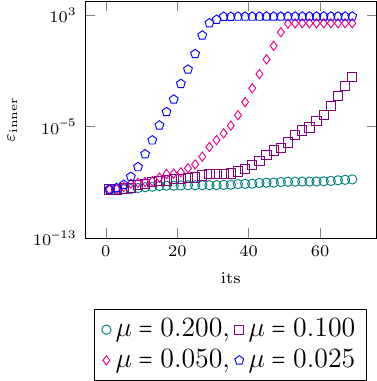}

		\caption{$\varepsilon=10^{-8}$}
		\label{fig:einner2}
	\end{subfigure}
	\hfill
	\begin{subfigure}[T]{0.3\columnwidth}
	
	\includegraphics{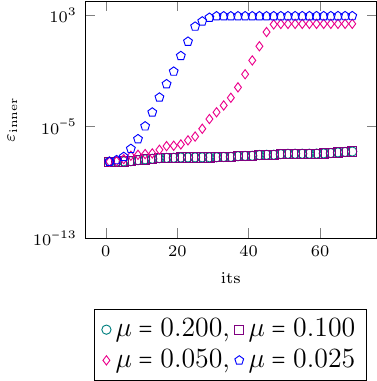}

		\caption{$\varepsilon=10^{-6}$}
		\label{fig:einner3}
	\end{subfigure}
	\caption{Value of $\varepsilon_{\text{inner}}^{(i)}$ as in \eqref{eq:pj-ineq} for different values of $\varepsilon$ when solving \eqref{eq:unshifted-trunc}, $m \rightarrow \infty$ with Algorithm~\ref{alg:FIGMRES}. Here $A(\mu)$ as in \eqref{eq:delay-A} with $n=1000$, random perturbations as in \eqref{eq:z-pert}, $\ell_j = 1$.}
	\label{fig:einner}
\end{figure}

\begin{figure}
\includegraphics{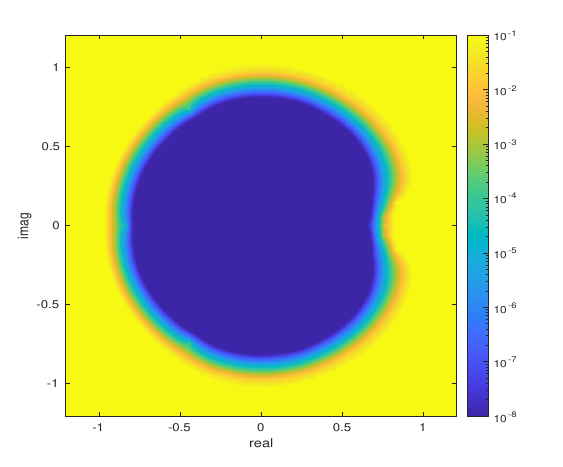}

\caption{Relative residual norm as in \eqref{eq:rel-res-rand} with $i=150$ when solving \eqref{eq:unshifted-trunc}, $m \rightarrow \infty$ with Algorithm~\ref{alg:FIGMRES} for a set of parameters $\{ \mu_1, \ldots, \mu_{k} \} \subset \mathbb{C}$. Here $A(\mu)$ as in \eqref{eq:delay-A} with $n=100$, random perturbations as in \eqref{eq:z-pert}, $\varepsilon_{\text{inner}}^{(i)}$ as in \eqref{eq:pj-ineq} and Algorithm~\ref{alg:FIGMRES} is executed once for $\mu_1 = 0.5 + 0.5i$, and we select $\ell_j = 1$, $\varepsilon=10^{-12}$.}
\label{fig:conv-radius}
\end{figure}

Figure~\ref{fig:asymfigs} illustrates the convergence of Algorithm~\ref{alg:FIGMRES} applied to \eqref{eq:unshifted-trunc} with $m \rightarrow \infty$, $A(\mu)$ as in \eqref{eq:delay-A}. Here the algorithm is run one time for the largest value of $\mu$ and then the function handle $\tilde{x}_i(\mu)$ is evaluated for each of the smaller values of $\mu$ considered. Each evaluation of $\tilde{x}_i(\mu)$ corresponds to solving one small linear least squares problem. As the algorithm is not run more than once, we generate a new approximation without modifying the random perturbations. Thus, for a given choice of $\mu$, the random perturbations introduced into Algorithm~\ref{alg:FIGMRES} will not disrupt convergence for smaller values of $\mu$. As discussed previously, we can expect fast convergence for small values of $\mu$ (see \eqref{eq:conv-bound}.\footnote{Note: though we consider only the first $n$ elements of the solution vector for the linear system as in \eqref{eq:unshifted-trunc} here, the same result applies.}) The convergence of Algorithm~\ref{alg:FIGMRES} is displayed alongside the convergence of applying preconditioned GMRES to this problem, i.e., considering an exact solve with $A(0)$ on each iteration.

Figure~\ref{fig:rtilde} illustrates the convergence of Algorithm~\ref{alg:FIGMRES} applied to the same problem described above, for different values of the parameter $\varepsilon$ as in \eqref{eq:pj-ineq}. We show the results for four values of $\mu$. Figure~\ref{fig:einner} shows the value of $\varepsilon_{\text{inner}}^{(i)}$ as in \eqref{eq:pj-ineq} for the same selection of $\varepsilon$ and $\mu$. In general, a smaller value for $\varepsilon$ leads to a lower relative residual norm at convergence. We can interpret this as, the more exact the application of the preconditioning step early in the algorithm, the smaller the relative residual norm is at convergence. It is clear that Algorithm~\ref{alg:FIGMRES} converges faster for smaller values of the parameter $\mu$. Figure~\ref{fig:einner} shows that if we execute our algorithm one time for the largest value of $\mu$, we can reuse the inner stopping criteria for smaller choices of $\mu$. This is because the inner stopping criteria becomes less restrictive for smaller values of $\mu$. Figures~\ref{fig:rtilde} and \ref{fig:einner} together show the inverse relationship between the quantities $\tilde{r}_i$ and $\varepsilon_{\text{inner}}^{(i)}$ for a fixed value of $\varepsilon$, as indicated in \eqref{eq:inverse}.

Figure~\ref{fig:conv-radius} displays the relative residual at convergence for the same problem as described above, for a set of parameters $\{ \mu_1, \ldots, \mu_k \} \subset \mathbb{C}$. Here Algorithm~\ref{alg:FIGMRES} is executed just once for $\mu_1 = 0.5 + 0.5i \in \mathbb{C}$. In this way, the random perturbations as in \eqref{eq:z-pert} correspond to the parameter $\mu_1$. We see that convergence of the approximation to \eqref{eq:our-prob} does not deteriorate for values $\mu_{l}$ such that $| \mu_{l} | \leq | \mu_1 |$, but does for $| \mu_{l} | > | \mu_1 |$. The degree of the deterioration in convergence depends on the distance of each particular $\mu_{l}$ from $\mu_1$.
Although our theory is for real $A(\mu)$ we have used complex $\mu$ in
Figure~\ref{fig:conv-radius} for visualization purposes.

\subsection{A finite element discretization of the parameterized Helmholtz equation} \label{sec:Helmholtz-sec}
\begin{figure}[h!]
	\begin{subfigure}[t]{0.45\columnwidth}
		\includegraphics[scale=.4]{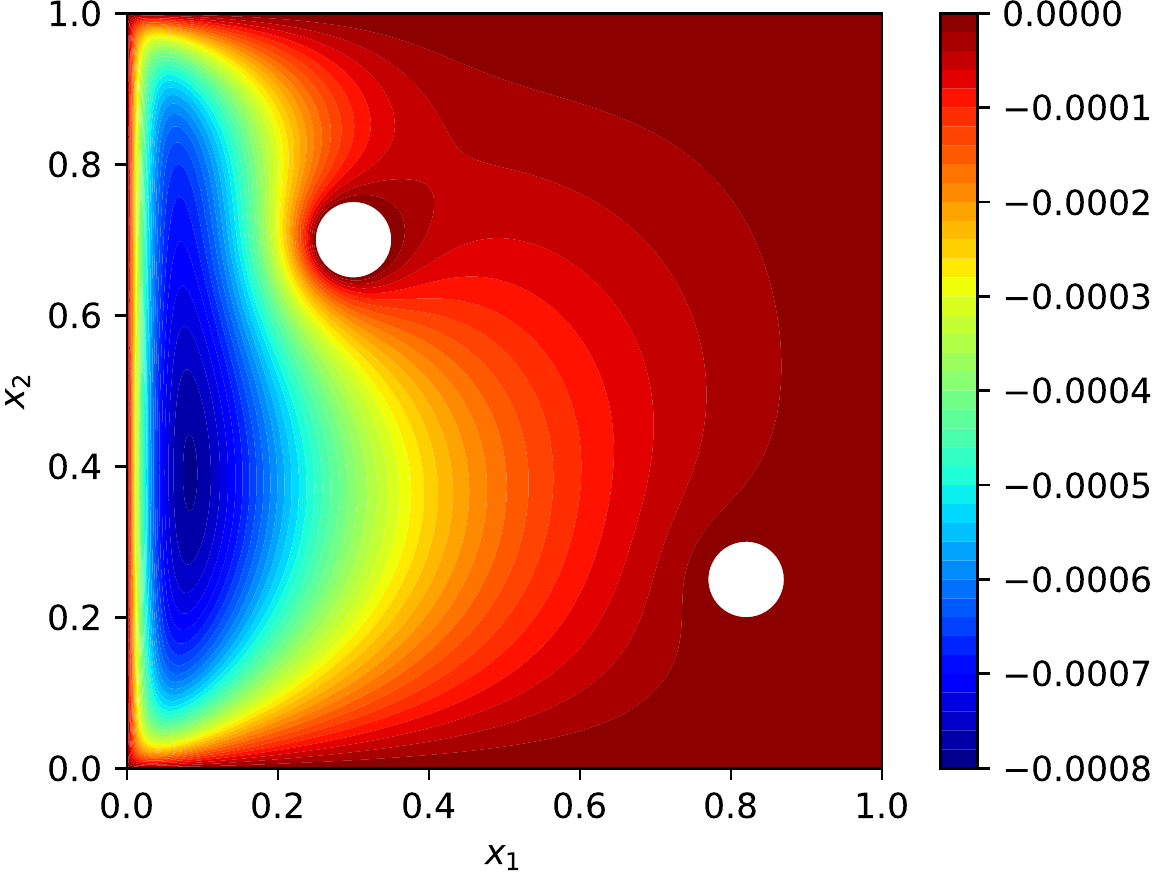}
		\caption{$\mu=.01$}
		\label{fig:numhelm1}
	\end{subfigure}
	\hfill
	\begin{subfigure}[t]{0.45\columnwidth}
		\includegraphics[scale=.4]{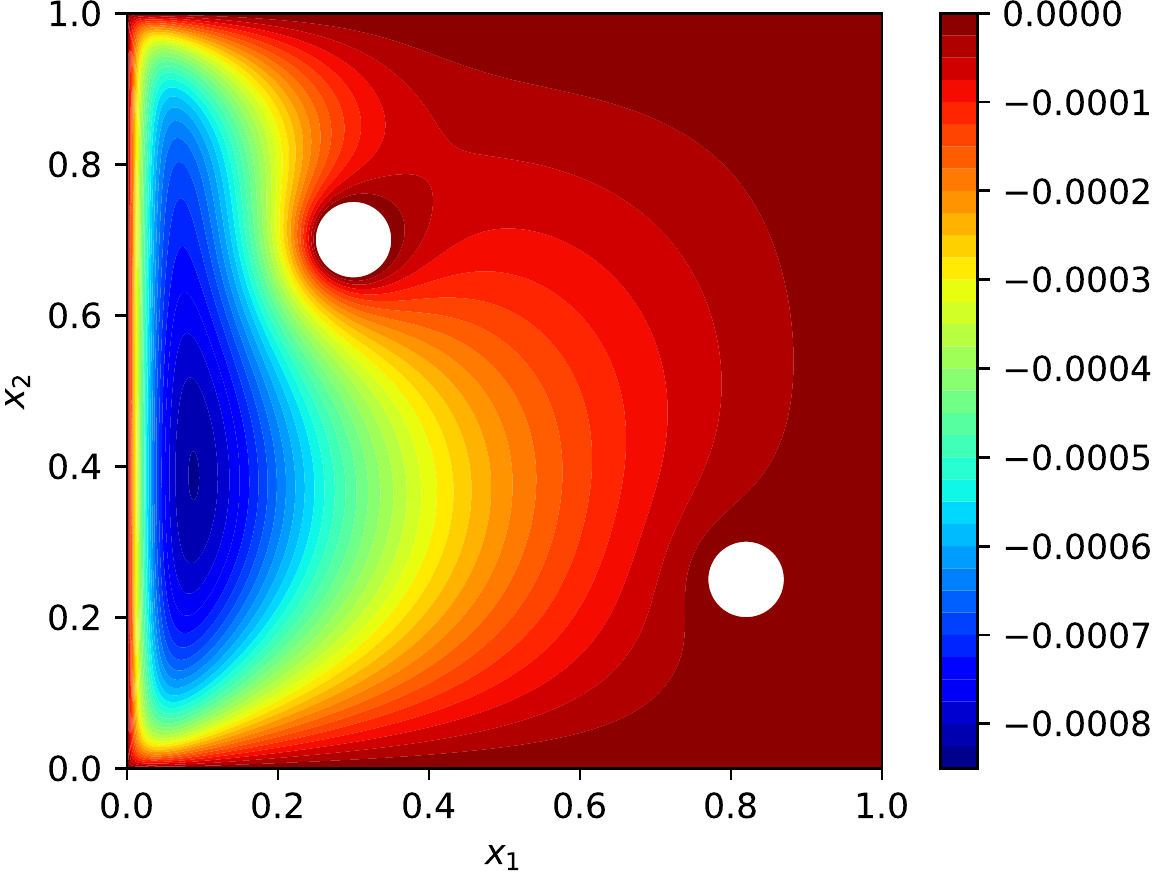}
		\caption{$\mu=1$}
		\label{fig:numhelm1a}
	\end{subfigure}

	\begin{subfigure}[t]{0.45\columnwidth}
		\includegraphics[scale=.4]{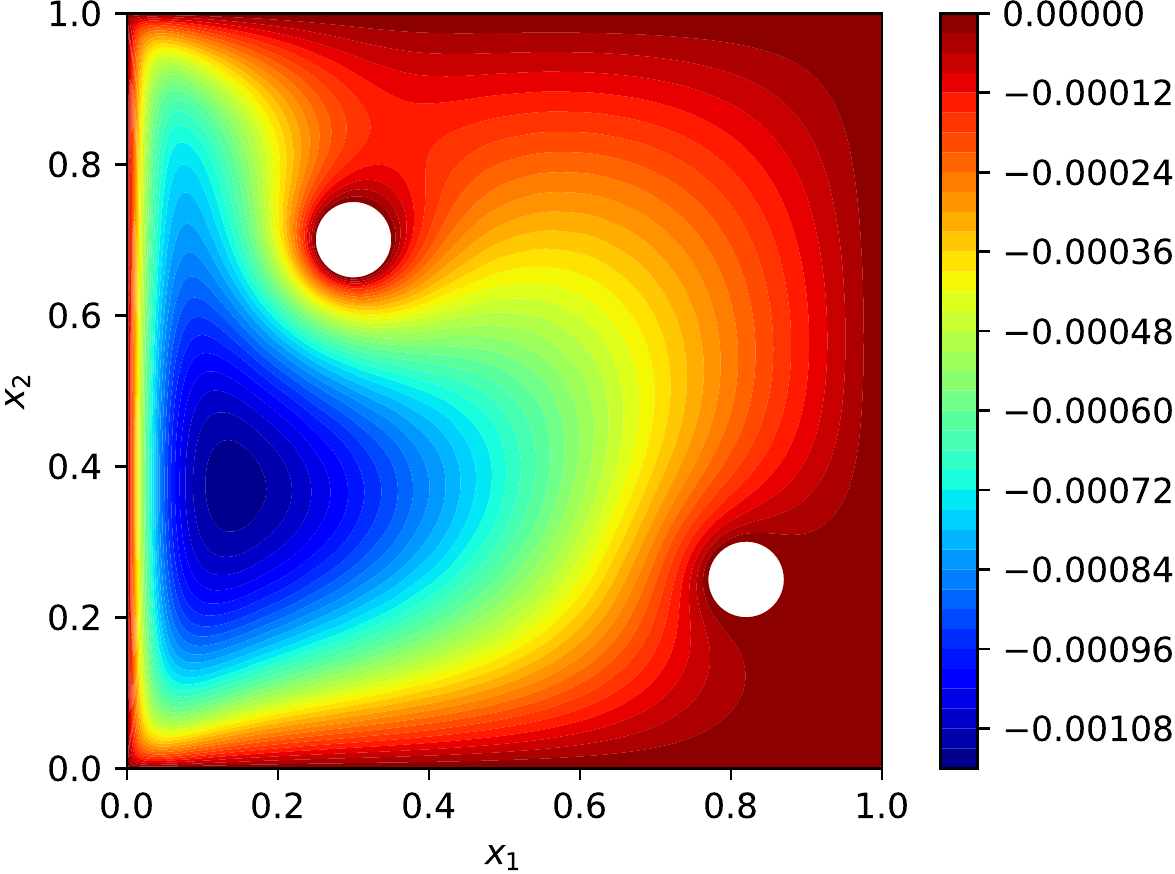}
		\caption{$\mu=2$}
		\label{fig:numhelm1b}
	\end{subfigure}
	\hfill
	\begin{subfigure}[t]{0.45\columnwidth}
		\includegraphics[scale=.4]{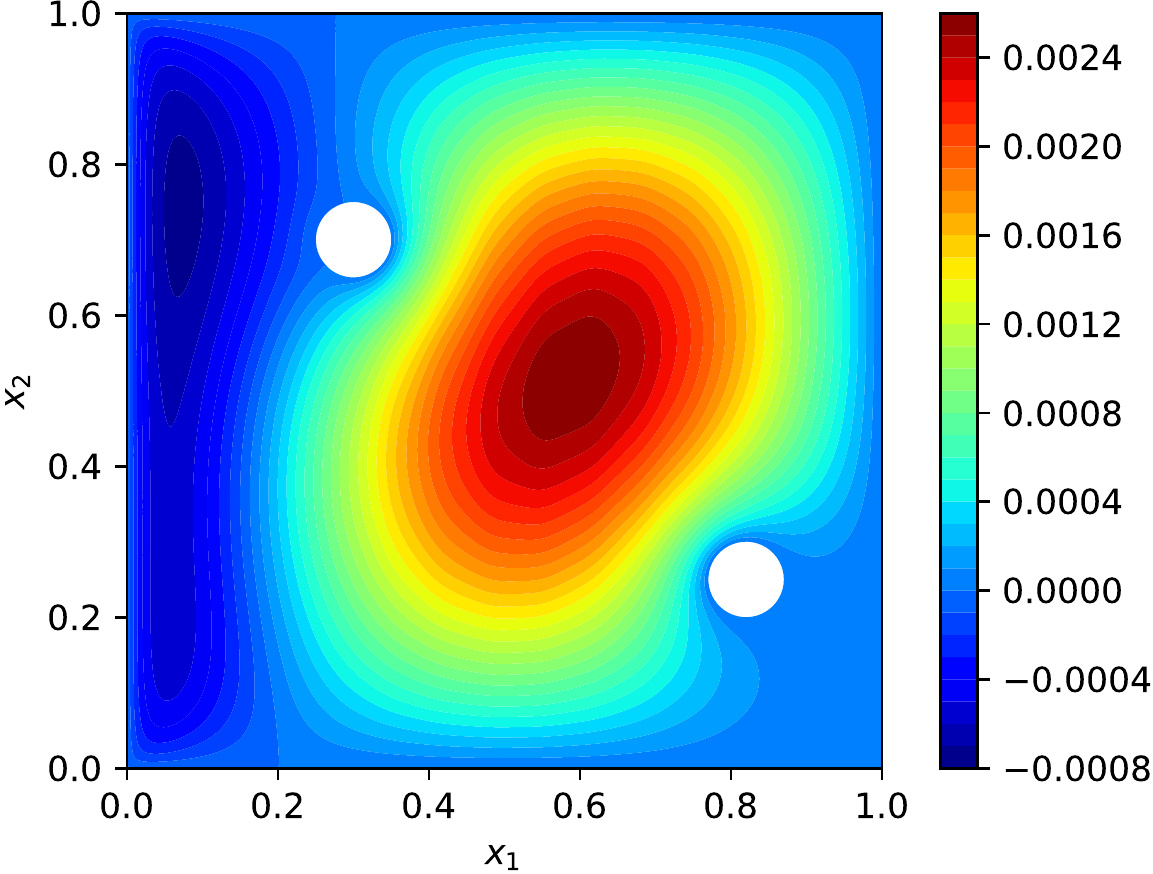}
		\caption{$\mu=2.5$}
		\label{fig:numhelm2}
	\end{subfigure}
	\caption{Numerical solutions to \eqref{eq:helmholtz_2d} generated by Algorithm~\ref{alg:FIGMRES}, $n = 245464$}
	\label{fig:Helm1}
\end{figure}

To illustrate the competitiveness of our proposed method, we consider approximating the solution $u(x)$ to the Helmholtz equation with a homogeneous Dirichlet boundary condition and a parameterized material coefficient given by
\begin{subequations} \label{eq:2dHelm}
\begin{alignat}{2}
   \left( \nabla^2 + f_1(\mu) \big( 1 + \mu k(x) \big) ^2 + f_2(\mu) \beta (x) \right) u(x) &= h(x), \qquad &&x \in \Omega, \\
   u(x) &= 0, &&x \in \partial \Omega,
\end{alignat}
\end{subequations}
where
\begin{gather*}
    k(x) =
    \begin{cases}
    1 + (x_1) \sin(\alpha \pi x_1), & x_1 \in [0, \frac{1}{2}), \\
    1 + (1 - x_1) \sin(\alpha \pi x_1), & x_1 \in [\frac{1}{2},1],
    \end{cases}
\end{gather*}
$h(x) = e^{-\alpha x_1}$, $\beta (x) = \sin(2 \pi x_1)$ for $x=(x_1,x_2)$, $f_1(\mu) = \mu$, $f_2(\mu) = \sin(\mu)$, $\alpha=30$, and the domain $\Omega$ is the unit square with two circular holes (see Figure~\ref{fig:Helm1}). The discretized problem is of the form
\begin{align}  \label{eq:helmholtz_2d}
	A_n(\mu) u_n(\mu) = h_n,
\end{align}
where
\begin{align} \label{eq:An-2d}
	A_n(\mu) := A_0 + \mu A_1 + 2 \mu^2 A_2 + \mu^3 A_3 + \sin (\mu) A_4
\end{align}
with $A_{\ell}$ large, sparse matrices arising from the finite element discretization of \eqref{eq:2dHelm} and $h_n$ the corresponding load vector.\footnote{The matrices and vector in the following simulations were generated using the finite element software FEniCS \cite{Alnaes:2015:FENICS}.} Thus, we use Algorithm~\ref{alg:FIGMRES}, returning a parametrized approximation to the solution in \eqref{eq:helmholtz_2d} which is cheap to evaluate for many values of $\mu$. Note that $A_n(\mu)$ as in \eqref{eq:An-2d} is nonlinear in $\mu$ and has an infinite Taylor series expansion.

To highlight the performance of our newly proposed method, we approximate the solution to the linear system \eqref{eq:unshifted-trunc} with $A_n(\mu)$ as in \eqref{eq:helmholtz_2d} and $m \rightarrow \infty$, in a preconditioned GMRES setting. In this way, we consider $A_n(\mu)$ without neglecting any of the terms in the infinite Taylor series expansion (see Section~\ref{sec:ext-to-infty}). In particular, we examine the convergence of three different iterative methods. The first two methods consider an inexact application of the preconditioning matrix $K_m^{-1}$ as in \eqref{eq:Km-inv} on iteration $i$. In both cases, the inexact application of the preconditioner is performed in a way which fulfills the dynamic inner stopping criteria as described in Section~\ref{sec:inner-tol}. The third method considers an exact application of the preconditioner. The methods are described as follows.

In the first inexact setting, we approximate the action of $A_0^{-1}$ in $K_m^{-1}$ with the iterative method Aggregation-based algebraic multigrid\footnote{Yvan Notay, AGMG software and documentation; see \texttt{http://agmg.eu}} (AGMG) \cite{NapovNotay2012,Notay2010,Notay2012} for a few (outer) iterations before applying the action of the identity matrix in place of $A_0^{-1}$. Specifically, we solve the linear systems with AGMG when the bound indicated in \eqref{eq:pj-ineq} is relatively restrictive, i.e. when $\varepsilon_{\text{inner}}^{(i)}$ is small. Once $\varepsilon_{\text{inner}}^{(i)}$ is sufficiently large, the preconditioner is instead applied via a very cheap approximation to $K_m^{-1}$ given by
\begin{align} \label{eq:last-precond}
	\begin{bmatrix}
	I & -A_1 & \cdots & -A_m  \\
	& I & & \\
	& & \ddots & \\
	& & & I
	\end{bmatrix} \in \mathbb{R}^{(m+1)n \times (m+1)n},
\end{align}
while still satisfying \eqref{eq:pj-ineq}. This option offers a clear computational advantage, removing a linear solve of size $n \times n$ on each iteration, while still retaining the assumptions of Lemma~\ref{matvec}, i.e., convergence of the inexact algorithm does not deteriorate as a result of this choice. The second inexact setting considers approximating the action of $A_0^{-1}$ with the method BiCGSTAB \cite{Vorst1992}. Specifically, here we specify the tolerance $tol$ in the iterative method with $tol = \varepsilon_{\text{inner}}^{(i-1)}$ on iteration $i$. In this way, the application of the preconditioner becomes more inexact with each iteration. This choice also satisfies the bound in \eqref{eq:pj-ineq} and convergence is achieved within a small number of iterations relative to the dimension $n$.

In the exact setting $A_0^{-1}$ is applied via an LU decomposition performed once before the execution of the algorithm, though other methods based on a decomposition of $A_0$ could also have been suitable here. Our aim is to compare the convergence of Algorithm~\ref{alg:FIGMRES} with dynamic inner stopping criteria to an exact preconditioned GMRES method when solving the discretized problem.

\begin{figure}[h!]
	\begin{subfigure}[t]{0.45\columnwidth}

	\includegraphics{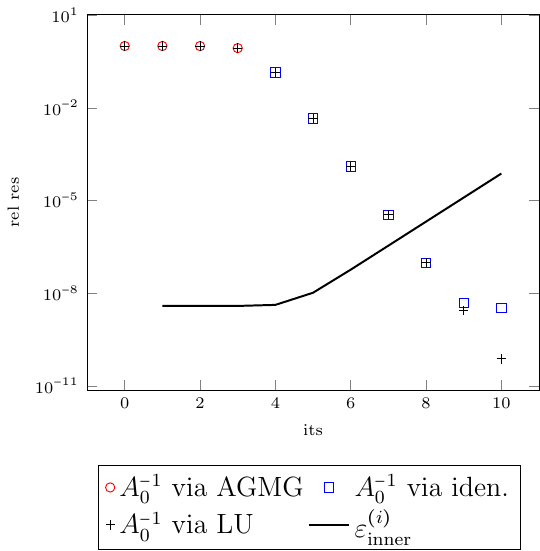}

    	\caption{$\mu=.01$, $s=.06$}
	\label{fig:alt_idea1a}
	\end{subfigure}
	\hfill
	\begin{subfigure}[t]{0.45\columnwidth}
	
	\includegraphics{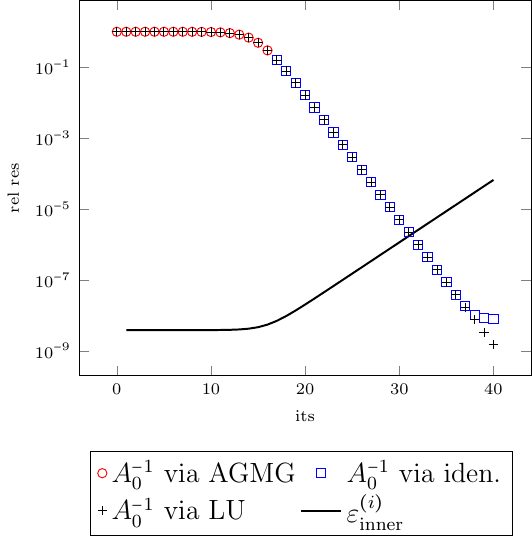}

		\caption{$\mu=1$, $s=1.5$}
		\label{fig:alt_idea1b}
	\end{subfigure}
\caption{Convergence for approximating \eqref{eq:helmholtz_2d} with Algorithm~\ref{alg:FIGMRES}. Here relative residual as in \eqref{eq:rel-res-3}, $\varepsilon=10^{-12}$, $n = 245464$ with $1713552$ nonzero elements in each $A_{\ell}$, $\ell = 0, . . . , 4$. The action of $A_0^{-1}$ in the preconditioner is approximated with the identity matrix whenever this substitute fulfills the inner stopping criteria \eqref{eq:pj-ineq}. The parameter $s$ is as described in Remark~\ref{scaling-factor}.}
\label{fig:alt_idea1}
\end{figure}

\begin{figure}
		\begin{subfigure}[t]{0.45\columnwidth}
	\includegraphics{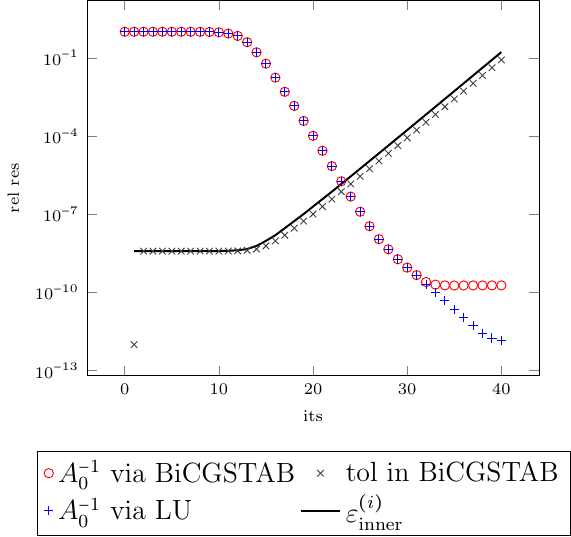}	

	\caption{$\mu=2$, $s=4$}
	\label{fig:alt_idea2a}
	\end{subfigure}
	\hfill
	\begin{subfigure}[t]{0.45\columnwidth}

	\includegraphics{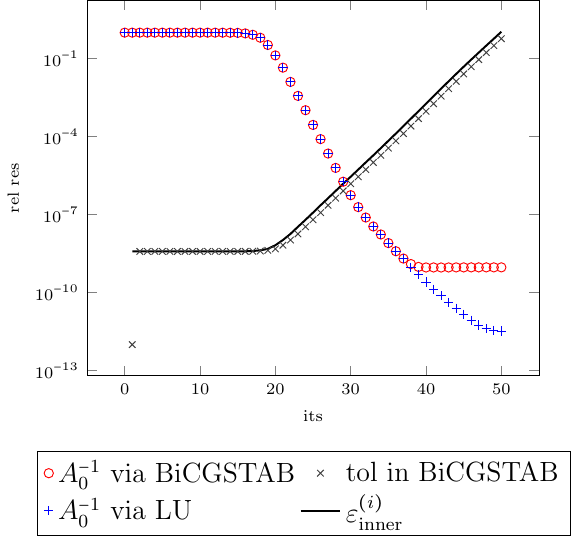}

	\caption{$\mu=2.5$, $s=5$}
	\label{fig:alt_idea2b}
	\end{subfigure}
\caption{Convergence for approximating \eqref{eq:helmholtz_2d} with Algorithm~\ref{alg:FIGMRES}. Here relative residual as in \eqref{eq:rel-res-3}, $\varepsilon=10^{-12}$, $n = 245464$ with $1713552$ nonzero elements in each $A_{\ell}$, $\ell = 0, . . . , 4$. The action of $A_0^{-1}$ in the preconditioner is approximated with the iterative method BiCGSTAB with tolerance given by $\varepsilon_{\text{inner}}^{(i-1)}$ on iteration $i$, fulfilling the inner stopping criteria \eqref{eq:pj-ineq}. The parameter $s$ is as described in Remark~\ref{scaling-factor}.}
\label{fig:alt_idea2}
\end{figure}

Figure~\ref{fig:Helm1} shows the numerical solutions generated by Algorithm~\ref{alg:FIGMRES} for four values of the parameter $\mu$. Figure~\ref{fig:alt_idea1} displays the convergence of the first inexact method as well as the exact method described above for two values of $\mu$, and Figure~\ref{fig:alt_idea2} shows an analogous figure for two values of $\mu$, where the second inexact method was used. Here we visualize the relative residual of the discretized problem
\begin{align} \label{eq:rel-res-3}
	\frac{\norm{A_n(\mu) \tilde{x}_i(\mu) - h_n}}{\norm{h_n}}
\end{align}
on iteration $i$, where $\tilde{x}_i(\mu)$ is the approximation given by Algorithm~\ref{alg:FIGMRES} or analogously, the approximation generated when applying the $K_m^{-1}$ exactly. The plots in Figure~\ref{fig:alt_idea1} are labeled to reflect when $K_m^{-1}$ is approximated by \eqref{eq:last-precond}, and the plots in Figure~\ref{fig:alt_idea2} show the tolerance $tol$ of the BiCGSTAB on each iteration. We note that the value of $tol$ increases inversely with the residual of the outer method.

To emphasize the advantage of using Algorithm~\ref{alg:FIGMRES} to approximate the solution to \eqref{eq:our-prob} for many values of $\mu$, consider the following example. Execute Algorithm~\ref{alg:FIGMRES} where $\mu=1$, $n=978452$, $s=5$ (see Remark~\ref{scaling-factor}), and apply the preconditioner $K_m^{-1}$ using the iterative method BiCGSTAB with growing tolerance on each iteration, analogous to the experiments in Figure~\ref{fig:alt_idea2}. The relative residual is approximately $6.8 \times 10^{-10}$ after $15$ outer iterations, taking approximately 2500 CPU seconds. Compare this to solving \eqref{eq:our-prob} for the particular value $\mu=1$, using the iterative method BiCGSTAB with tolerance set to $10^{-10}$. This process takes approximately 250 CPU seconds, but returns an approximation for just one value of the parameter. Recall that Algorithm~\ref{alg:FIGMRES} requires fewer outer iterations to achieve convergence for smaller values of $\mu$. Thus, $15$ outer iterations of Algorithm~\ref{alg:FIGMRES} gives access to the approximation to \eqref{eq:our-prob} for values $|\mu_i| \leq 1$, $\mu_i \in \mathbb{R}$. Our novel method is more efficient when we require the solution to $\eqref{eq:our-prob}$ for more than $10$ such $\mu_i$. Analogous results can be obtained for different values of the parameter $\mu$.

\begin{remark}[Complexity of the inexact algorithm]
In the experiments visualized in Figures~\ref{fig:alt_idea1} and \ref{fig:alt_idea2}, Algorithm~\ref{alg:FIGMRES} offers a dramatic reduction in complexity over the exact preconditioned algorithm with the following setup. In the exact preconditioning setting we compute the LU decomposition of $A_0$ followed by two triangular solves on each outer iteration $i=1,\ldots,j$, i.e., for dense matrices $\mathcal{O} (2/3 n^3 + jn^2)$ operations. The inexact methods apply the action of $A_0^{-1}$ approximately, either via AGMG with complexity $\mathcal{O} (n)$ on each of the first few outer iterations or BiCGSTAB with complexity $\mathcal{O} (n)$ per iteration and growing tolerance. The three methods require the same matrix-vector multiplication, $\mathcal{O}(j^2n)$ and orthogonalization, $\mathcal{O}(j^3 n)$.
\end{remark}
\begin{remark}[Comparison to infinite GMRES]
Algorithm~\ref{alg:FIGMRES} requires approximately the same number of outer iterations as infinite GMRES \cite{JarlebringCorrenty1}, when controlling for the magnitude of $\mu$ and the dimension. Therefore, the inexact algorithm offers a dramatic improvement in terms of performace, as the decomposition of $A_0$ can be skipped entirely.
\end{remark}

\begin{remark}[Scaling factor]  \label{scaling-factor}
Solving \eqref{eq:our-prob} with Algorithm~\ref{alg:FIGMRES} can be formulated equivalently as follows. Define
\[
	\hat{A}(\mu) := \sum_{\ell=0}^{\infty} A_{\ell} (s\mu)^{\ell},
\]
where $s$ is a constant such that $s>\mu$ and $A_{\ell}$ as in \eqref{eq:taylor-coeff}. Perform the linearization on $\hat{A}(\mu)$, noting that $\hat{A}(\mu/s)= A(\mu)$ and equivalently,
\[
	\hat{A}\left(\frac{\mu}{s}\right) \hat{x}\left(\frac{\mu}{s}\right) = b \iff A(\mu) x(\mu) = b.
\]
Hence, in this case Algorithm~\ref{alg:FIGMRES} returns $\hat{x}$ and we can easily recover the desired solution. In practice, this setting converges in fewer outer iterations.
\end{remark}

\section{Conclusions and outlook}
The main contribution of this paper is a new Krylov subspace method. This method solves parameterized systems of the form \eqref{eq:our-prob}, by returning a parameterized approximation $\tilde{x}$ which is cheap to evaluate for many values of $\mu$. A companion linearization was formed where $\mu$ appeared only linearly, and one orthonormal basis matrix for the Krylov subspace \eqref{eq:our-krylov} was built in way which included all the terms in the infinite Taylor series expansion of $A(\mu)$. Analysis and numerical experiments showed that the convergence of the method can be predicted by the magnitude of the parameter $\mu$.

The new GMRES method incorporates preconditioning where the preconditioning matrix is of the form \eqref{eq:Km-inv}. In applying this preconditioner, a linear solve of dimension $n \times n$ is required on each iteration. Our framework considers a flexible setting for the preconditioning, allowing for an inexact solve determined by a dynamic inner stopping criteria. This inexactness offers a substantial reduction in complexity compared to the previously proposed method infinite GMRES \cite{JarlebringCorrenty1}. Numerical simulations verify that this modification does not degrade convergence. Additionally, the method is suitable for large-scale simulations.

This work can be extended in various ways beyond the scope of this paper. There are several variants of the infinite Arnoldi method, e.g. the Chebyshev version \cite{Jarlebring:2012:INFARNOLDI} and restarting variations \cite{Jarlebring:2014:SCHUR}. These strategies could be applied to the methods presented in this paper, but this would require adaptation based on the specific structure of the problem and therefore, further analysis.


\bibliographystyle{plain}
\bibliography{siobhanbib,eliasbib}

\begin{thebibliography}{10}

\bibitem{AhujaEtAl2015}
K.~Ahuja, P.~Benner, E.~de~Sturler, and L.~Feng.
\newblock Recycling {B}i{CGSTAB} with an application to parametric model order
  reduction.
\newblock {\em SIAM J. Sci. Comput.}, 37(5):S429--S446, 2015.

\bibitem{AhujaEtAl2012}
K.~Ahuja, E.~de~Sturler, S.~Gugercin, and E.~R. Chang.
\newblock Recycling {B}i{CG} with an application to model reduction.
\newblock {\em SIAM J. Sci. Comput.}, 34(4):A1925--A1949, 2012.

\bibitem{Alnaes:2015:FENICS}
M.~S. Alnaes, J.~Blechta, J.~Hake, A.~Johansson, B.~Kehlet, A.~Logg,
  C.~Richardson, J.~Ring, M.~E. Rognes, and G.~N. Wells.
\newblock The {FEniCS} project version 1.5.
\newblock {\em Archive of Numerical Software}, 3, 2015.

\bibitem{AntoulasEtAl2020}
A.C. Antoulas, C.A. Beattie, and S.~Gugercin.
\newblock {\em Interpolatory methods for model reduction}.
\newblock SIAM, 2020.

\bibitem{Axelsson1991}
O.~Axelsson and P.S. Vassilevski.
\newblock A black box generalized conjugate gradient solver with inner
  iterations and variable-step preconditioning.
\newblock {\em SIAM J. Matrix Anal. Appl.}, 12:625--644, 1991.

\bibitem{BaglamaCal1998}
J.~Baglama, D.~Valvetti, G.~H. Golub, and L.~Reichel.
\newblock Adaptively preconditioned {GMRES} algorithms.
\newblock {\em SIAM J. Sci. Comput.}, 20:243--269, 1998.

\bibitem{Bahkos:2017:MULTIPREC}
T.~Bakhos, P.~K. Kitanidis, S.~Ladenheim, A.~K. Saibaba, and D.~B. Szyld.
\newblock Multipreconditioned {GMRES} for shifted systems.
\newblock {\em SIAM J. Sci. Comput.}, pages S222--S247, 2017.

\bibitem{Baumann2015NestedKM}
M.~Baumann and M.~B. van Gijzen.
\newblock Nested {K}rylov methods for shifted linear systems.
\newblock {\em SIAM J. Sci. Comput.}, 37, 2015.

\bibitem{BaylissEtAl1983}
A.~Bayliss, C.I. Goldstein, and E.~Turkel.
\newblock An iterative method for the {H}elmholtz equation.
\newblock {\em J. Comput. Phys.}, 49:443--457, 1983.

\bibitem{BeattieGug09}
C.~Beattie and S.~Gugercin.
\newblock Interpolatory projection methods for structure-preserving model
  reduction.
\newblock {\em Syst. Control Lett.}, 58(3):225--232, 2009.

\bibitem{Bouras2000ARS}
A.~Bouras and V.~Frayss{\'e}.
\newblock A relaxation strategy for inexact matrix-vector products for {K}rylov
  methods.
\newblock Technical Report TR/PA/00/15, CERFACS, Toulouse, France, 2000.

\bibitem{Bouras00arelaxation}
A.~Bouras, V.~Fraysse, and L.~Giraud.
\newblock A relaxation strategy for inner-outer linear solvers in domain
  decomposition methods.
\newblock Technical Report TR/PA/00/16, CERFACS, Toulouse, France, 2000.

\bibitem{Bouras00arelaxation2}
A.~Bouras and V.~Frayssé.
\newblock A relaxation strategy for the {A}rnoldi method in eigenproblems.
\newblock Technical report TR/PA/00/17, CERFACS, Toulouse, France, 2000.

\bibitem{Brown91}
P.~N. Brown.
\newblock A theoretical comparison of the {A}rnoldi and {GMRES} algorithms.
\newblock {\em SIAM J. Sci. Stat. Comput.}, 12(1):58--78, 1991.

\bibitem{Chapman1998}
A.~Chapman and Y.~Saad.
\newblock Deflated and augmented {K}rylov subspace techniques.
\newblock {\em Numer. Linear Algebra Appl.}, 4:43--66, 1998.

\bibitem{deSturler1996}
E.~de~Sturler.
\newblock Nested {K}rylov methods based on {GCR}.
\newblock {\em J. Comput. Appl. Math.}, 67:15--41, 1996.

\bibitem{SturlerEtAl15}
E.~de~Sturler, S.~Gugercin, M.~E. Kilmer, S.~Chaturantabut, C.~Beattie, and
  M.~O'Connell.
\newblock Nonlinear parametric inversion using interpolatory model reduction.
\newblock {\em SIAM J. Sci. Comput.}, 37(3):B495--B517, 2015.

\bibitem{ElmanErn2001}
H.~C. Elman, O.~G. Ernst, and D.~P. O'Leary.
\newblock A multigrid method enhanced by {K}rylov subspace iteration for
  discrete {H}elmholtz equations.
\newblock {\em SIAM J. Sci. Comput.}, 23:1291--1315, 2001.

\bibitem{ElmanEtAl2001}
H.C. Elman, O.G. Ernst, and D.P. O'Leary.
\newblock A multigrid method enhanced by {K}rylov subspace iteration for
  discrete {H}elmholtz equations.
\newblock {\em SIAM J. Sci. Comput.}, 23(4):1291--1315, 2001.

\bibitem{Erhel1996}
J.~Erhel, K.~Buggage, and B.~Pohl.
\newblock Restarted {GMRES} preconditioned by deflation.
\newblock {\em J. Comput. Appl. Math.}, 69:303--318, 1996.

\bibitem{Erlangga08}
Y.A. Erlangga.
\newblock Advances in iterative methods and preconditioners for the {H}elmholtz
  equation.
\newblock {\em Archives Comput. Methods Engin.}, 15:37--66, 2008.

\bibitem{ErlanggaEtAl2004}
Y.A. Erlangga, C.~Vuik, and C.W. Oosterlee.
\newblock On a class of preconditioners for solving the {H}elmholtz equation.
\newblock {\em Appl. Numer. Math.}, 50:409--425, 2004.

\bibitem{FengEtAl2013}
L.~Feng, P.~Benner, and J.~G. Korvink.
\newblock Subspace recycling accelerates the parametric macro-modeling of
  {MEMS}.
\newblock {\em Int. J. Numer. Methods Eng.}, 94(1):84--110, 2013.

\bibitem{ShiftedFreund}
R.~W. Freund.
\newblock {\em Solution of shifted linear systems by quasi-minimal residual
  iterations}.
\newblock De Gruyter, 1993.

\bibitem{FrommerGlassner98}
A.~Frommer and U.~Glässner.
\newblock Restarted {GMRES} for shifted linear systems.
\newblock {\em SIAM J. Sci. Comput.}, 19(1):15--26, 1998.

\bibitem{FrommerMaass99:20}
A.~Frommer and P.~Maass.
\newblock Fast {CG}-based methods for {T}ikhonov--{P}hillips regularization.
\newblock {\em SIAM J. Sci. Comput.}, 20(5):1831--1850, 1999.

\bibitem{Golub1999InexactPC}
G.~Golub and Q.~Ye.
\newblock Inexact preconditioned conjugate gradient method with inner-outer
  iteration.
\newblock {\em SIAM J. Sci. Comput.}, 21:1305--1320, 1999.

\bibitem{GuSimoncini}
G.-D. Gu and V.~Simoncini.
\newblock Numerical solution of parameter-dependent linear systems.
\newblock {\em Numer. Linear Algebra Appl.}, 12(9):923--940, 2005.

\bibitem{Gu:2003:STABILITY}
K.~Gu, V.~Kharitonov, and J.~Chen.
\newblock {\em Stability of {T}ime-{D}elay {S}ystems}.
\newblock Control Engineering. Boston, MA: Birkh\"auser, 2003.

\bibitem{JarlebringCorrenty1}
E.~Jarlebring and S.~Correnty.
\newblock Infinite {GMRES} for parameterized linear systems.
\newblock Technical report, 2022.
\newblock Accepted for publication in \textit{SIAM J. Matrix Anal. Appl.}

\bibitem{Jarlebring:2014:SCHUR}
E.~Jarlebring, K.~Meerbergen, and W.~Michiels.
\newblock Computing a partial {Schur} factorization of nonlinear eigenvalue
  problems using the infinite {Arnoldi} method.
\newblock {\em SIAM J. Matrix Anal. Appl.}, 35(2):411--436, 2014.

\bibitem{DelayProblem}
E.~Jarlebring, K.~Meerbergen, and W-Michiels.
\newblock An {A}rnoldi method with structured starting vectors for the delay
  eigenvalue problem.
\newblock {\em Proceedings of the 9th IFAC Workshop on Time Delay Systems,
  Prague}, 2010.

\bibitem{TensorArnoldi}
E.~Jarlebring, G.~Mele, and O.~Runborg.
\newblock The waveguide eigenvalue problem and the tensor infinite {A}rnoldi
  method.
\newblock {\em SIAM J. Sci. Comput.}, 39(3), 2017.

\bibitem{Jarlebring:2012:INFARNOLDI}
E.~Jarlebring, W.~Michiels, and K.~Meerbergen.
\newblock A linear eigenvalue algorithm for the nonlinear eigenvalue problem.
\newblock {\em Numer. Math.}, 122(1):169--195, 2012.

\bibitem{Joubert1994}
W.~Joubert.
\newblock A robust {GMRES}-based adaptive polynomial preconditioning algorithm
  for nonsymmetric linear systems.
\newblock {\em SIAM J. Sci. Comput.}, 15:427--439, 1994.

\bibitem{KilmerOleary01}
M.~E. Kilmer and D.~P. O'Leary.
\newblock Choosing regularization parameters in iterative methods for ill-posed
  problems.
\newblock {\em SIAM J. Matrix Anal. Appl.}, 22(4):1204--1221, 2001.

\bibitem{KilmerSturler06}
M.~E. Kilmer and E.~Sturler.
\newblock Recycling subspace information for diffuse optical tomography.
\newblock {\em SIAM J. Sci. Comput.}, 27(6):2140--2166, 2006.

\bibitem{Kressner:2014:TOAR}
D.~Kressner and J.~Roman.
\newblock Memory-efficient {Arnoldi} algorithms for linearizations of matrix
  polynomials in {Chebyshev} basis.
\newblock {\em Numer. Linear Algebra Appl.}, 21(4):569--588, 2014.

\bibitem{KressnerToblerLR}
D.~Kressner and C.~Tobler.
\newblock Low-rank tensor {K}rylov subspace methods for parametrized linear
  systems.
\newblock {\em SIAM J. Matrix Anal. Appl.}, 32:1288--1316, 2011.

\bibitem{Mackey:2006:VECT}
S.~Mackey, N.~Mackey, C.~Mehl, and V.~Mehrmann.
\newblock Vector spaces of linearizations for matrix polynomials.
\newblock {\em SIAM J. Matrix Anal. Appl.}, 28:971--1004, 2006.

\bibitem{Michiels2011}
W.~Michiels, E.~Jarlebring, and K.~Meerbergen.
\newblock Krylov-based model order reduction of time-delay systems.
\newblock {\em SIAM J. Matrix Anal. Appl.}, 32(4):1399--1421, 2011.

\bibitem{Michiels:2007:STABILITYBOOK}
W.~Michiels and S.-I. Niculescu.
\newblock {\em Stability and Stabilization of Time-Delay Systems: An
  Eigenvalue-Based Approach}.
\newblock Advances in Design and Control 12. SIAM Publications, Philadelphia,
  2007.

\bibitem{NapovNotay2012}
A.~Napov and Y.~Notay.
\newblock An algebraic multigrid method with guaranteed convergence rate.
\newblock {\em SIAM J. Sci. Comput.}, 34:A1079--A1109, 2012.

\bibitem{Notay2000}
Y.~Notay.
\newblock Flexible conjugate gradient.
\newblock {\em SIAM J. Sci. Comput.}, 22:1444--1460, 2000.

\bibitem{Notay2010}
Y.~Notay.
\newblock An aggregation-based algebraic multigrid method.
\newblock {\em Electronic Transactions on Numerical Analysis}, 37:123--146,
  2010.

\bibitem{Notay2012}
Y.~Notay.
\newblock Aggregation-based algebraic multigrid for convection-diffusion
  equations.
\newblock {\em SIAM J. Sci. Comput.}, 34:A2288--A2316, 2012.

\bibitem{ParksEtAl2006}
M.~L. Parks, E.~de~Sturler, G.~Mackey, D.~D. Johnson, and S.~Maiti.
\newblock Recycling {K}rylov subspaces for sequences of linear systems.
\newblock {\em SIAM J. Sci. Comput.}, 28(5):1651--1674, 2006.

\bibitem{Saad1993}
Y.~Saad.
\newblock A flexible inner-outer preconditioned {GMRES} algorithm.
\newblock {\em SIAM J. Sci. Comput.}, 14(2):461--469, 1993.

\bibitem{Saad:1986:GMRES}
Y.~Saad and M.~H. Schultz.
\newblock {GMRES}: A generalized minimal residual algorithm for solving
  nonsymmetric linear systems.
\newblock {\em SIAM J. Sci. Stat. Comput.}, 7:856--869, 1986.

\bibitem{Simoncini2002a}
V.~Simoncini and L.~Eld\'en.
\newblock Inexact {Rayleigh} quotient-type methods for eigenvalue computations.
\newblock {\em BIT}, 42(1):159--182, 2002.

\bibitem{inexKry}
V.~Simoncini and D.~Szyld.
\newblock Theory of inexact {K}rylov subspace methods and applications to
  scientific computing.
\newblock {\em SIAM J. Sci. Comput.}, 25, May 2002.

\bibitem{SimonciniSzyld2003}
V.~Simoncini and D.~B. Szyld.
\newblock Flexible inner-outer {K}rylov subspace methods.
\newblock {\em SIAM J. Numer. Anal.}, 40:2219--2239, 2003.

\bibitem{etna_vol45_pp499-523}
K.~M. Soodhalter.
\newblock Two recursive {GMRES}-type methods for shifted linear systems with
  general preconditioning.
\newblock {\em Electron. Trans. Numer. Anal.}, 45:499--523, 2016.

\bibitem{SOODHALTER2014105}
K.~M. Soodhalter, D.~B. Szyld, and F.~Xue.
\newblock Krylov subspace recycling for sequences of shifted linear systems.
\newblock {\em Appl. Numer. Math.}, 81:105--118, 2014.

\bibitem{Szyld2001}
D.~B. Szyld and J.~A. Vogel.
\newblock {FQMR}: A flexible quasi-minimal residual method with inexact
  preconditioning.
\newblock {\em SIAM J. Sci. Comput.}, 23:363--380, 2001.

\bibitem{VanBeeumen:2015:CORK}
R.~{V}an Beeumen, K.~Meerbergen, and W.~Michiels.
\newblock Compact rational {Krylov} methods for nonlinear eigenvalue problems.
\newblock {\em SIAM J. Sci. Comput.}, 36(2):820--838, 2015.

\bibitem{Vorst1992}
H.~A. van~der Vorst.
\newblock {Bi-CGSTAB}: A fast and smoothly converging variant of {Bi-CG} for
  the solution of nonsymmetric linear systems.
\newblock {\em SIAM J. Sci. Comput.}, 13(2):631--644, 1992.

\bibitem{vanDerVorstVuik1994}
H.~A. {V}an~der {V}orst and C.~Vuik.
\newblock {GMRESR}: A family of nested {GMRES} methods.
\newblock {\em Numer. Linear Algebra Appl.}, 1:369--386, 1994.

\bibitem{ZaranVarga}
R.~S. Varga.
\newblock A comparison of the successive overrelaxation method of
  semi-iterative methods using {C}hebyshev polynomials.
\newblock {\em J. Soc. Indust. Appl. Math.}, 5, 1957.

\bibitem{Vuik19932}
C.~Vuik.
\newblock Further experiences with {GMRESR}.
\newblock {\em Supercomputer}, 55:13--27, 1993.

\bibitem{vuik1993}
C.~Vuik.
\newblock New insights in {GMRES}-like methods with variable preconditioners.
\newblock {\em J. Comput. Appl. Math.}, 61:189--204, 1995.

\bibitem{WarsaBenzi}
J.~S. Warsa, M.~Benzi, T.~A. Wareing, and J.~E. Morel.
\newblock Preconditioning for a mixed discontinuous finite element method for
  radiation diffusion.
\newblock {\em Numer. Linear Algebra Appl.}, 11:795--811, 2004.

\end{thebibliography}

\end{document}